\documentclass[preprint,12pt]{elsarticle}

\usepackage{amssymb}
\usepackage{graphicx}
\usepackage{amsmath}
\usepackage{latexsym}
\usepackage{amssymb}
\usepackage{pstricks}
\usepackage{graphics,color}
\usepackage{color}
\usepackage{amsfonts}
\usepackage{subfigure}
\usepackage{subfigmat}

\newcommand{\newtext}[1]{#1}

\newcommand{\om}{\Omega}

\newcommand{\bfx}{\mbox{\boldmath{$x$}}}
\newcommand{\bfn}{\mbox{\boldmath{$n$}}}
\newcommand{\bff}{\mbox{\boldmath{$f$}}}

\newcommand{\bfG}{\mbox{\boldmath{$G$}}}

\newcommand{\bfh}{\mbox{\boldmath{$h$}}}

\newcommand{\bfJ}{\mbox{\boldmath{$J$}}}
\newcommand{\bfS}{\mbox{\boldmath{$S$}}}

\newcommand{\bfphi}{\mbox{\boldmath{$\phi$}}}

\newcommand{\bfxi}{\mbox{\boldmath{$\xi$}}}

\newcommand{\M}[1]{\boldsymbol{\mathsf{#1}}}
\newcommand{\de}{\,{\mathrm d}}
\newcommand{\trn}{^\mathsf{T}}
\DeclareMathOperator*{\assembly}{{\sf A}}

\newtheorem{thm}{Theorem}
\newtheorem{lem}[thm]{Lemma}
\newtheorem{dfn}[thm]{Definition}

\newdefinition{rmk}{Remark}
\newdefinition{crl}{Corollary}
\newproof{pf}{Proof}

\begin{document}

\begin{frontmatter}



\title{Analysis of coupled transport phenomena in concrete \\
at elevated temperatures}


\author[CIDEAS,KM]{\corref{cor1}Michal Bene\v{s}}
\cortext[cor1]{Corresponding address: Department of Mathematics,
Faculty of Civil Engineering, Czech Technical University in Prague,
Th\'{a}kurova 7, 166 29 Prague 6, Czech Republic}
\ead{xbenesm3@fsv.cvut.cz}

\author[KBK]{Radek \v{S}tefan}
\author[FCEKSM]{Jan Zeman}

\address[CIDEAS]{Centre for Integrated Design of Advanced Structures,}
\address[KM]{Department of Mathematics,}
\address[KBK]{Department of Concrete and Masonry Structures,}
\address[FCEKSM]{Department of Mechanics,\\Faculty of Civil Engineering,\\
Czech Technical University in Prague,\\
Th\'{a}kurova 7, 166 29 Prague 6, Czech Republic}

\begin{abstract}
In this paper, we study a non-linear numerical scheme arising from the implicit
time discretization of the Ba\v{z}ant-Thonguthai model for hygro-thermal
behavior of concrete at high temperatures. Existence and uniqueness of the
time-discrete solution in two dimensions is established using the theory of
pseudomonotone operators in Banach spaces. Next, the spatial discretization is
accomplished by the conforming finite element method. An illustrative numerical
example shows that the model reproduces well the rapid increase of
pore pressure in wet concrete due to extreme heating. Such phenomenon is of
particular interest for the safety assessment of concrete structures prone to
thermally-induced spalling.
\end{abstract}

\begin{keyword}
heat and \newtext{moisture} transfer phenomena in concrete at high
temperatures \sep Rothe method \sep pseudo-monotone operators \sep
finite element discretization \sep spalling
\end{keyword}

\end{frontmatter}

\section{Introduction}\label{Intro}

The hygro-thermal behavior of concrete plays a crucial role in the
assessment of the reliability and lifetime of concrete structures.
The heat and mass transfer processes become particularly important
at high temperatures, where the increased pressure in pores may lead
to catastrophic service failures. Since high-temperature experiments
are very expensive, predictive modeling of humidity migration and
pore pressure developments can result in significant economic
savings. The first mathematical models of concrete exposed to
temperatures exceeding $100^{\circ}$C were formulated by Ba\v{z}ant
and~Thonguthai in \cite{BaTh1978}. Since then, a considerable effort
has been invested into detailed numerical simulations of concrete
structures subject to high temperatures. However, much less
attention has been given to the qualitative properties of the model,
as well as of the related numerical methods.

In particular, the only related work the authors are aware of is due to
Dal\'{i}k \emph{et al.}~\cite{DaDaVa2000}, who analyzed the numerical solution
of the Kiessl model~\cite{Kiessl1983} for moisture and heat transfer in porous
materials. They proved some existence and regularity results and suggested an
efficient numerical approach to the solution of the resulting system of highly
non-liner equations. However, the Kiessl model is valid for limited temperature
range only and as such it is inappropriate for high-temperature applications. In
this contribution, we extend the work~\cite{DaDaVa2000} by proving the existence
and uniqueness of an approximate solution for the Ba\v{z}ant-Thonguthai model,
arising from the semi-implicit discretization in time. A fully discrete
algorithm is then obtained by standard finite element discretization and its
performance is illustrated for a model problem of a concrete segment exposed to
transient heating according to the standard ISO fire curve. Here, the focus is
on the short-term pore pressure build up, which is decisive for the assessment
of so-called thermal spalling during fire.

\newtext{%
At this point, it is fair to mention that the Ba\v{z}ant-Thonguthai
model was later extended towards more detailed multi-phase description, see e.g.
the works of Gawin \emph{et al.}~\cite{GaMaSch1999}, Tenchev \emph{et
al.}~\cite{TeLiPu2001} and Davie
\emph{et al.}~\cite{DaPeBi2006} for specific examples. 
When compared to the original version, these advanced models provide better
insight into physical and chemical processes in concrete (such as influence of
gel water, pore water, capillary water, chemical reactions at elevated
temperatures, etc.). Such potential increase in accuracy, however, comes at the
expense of increased number of parameters, which
typically reflect complex multi-scale nature of concrete. Hence, their
experimental determination is rather complicated and the parameters can often
only be calibrated by sub-scale simulations. Therefore, in this work we adopt a
pragmatic approach and consider the single-phase Ba\v{z}ant-Thonguthai model
with parameters provided by heuristic relations, obtained from regression of
reliable macroscopic experiments.} 

The paper is organized as follows. In Section \ref{model}, we
present the general single-phase\newtext{, purely macroscopic,} model for
prediction of hygro-thermal behavior of heated concrete. In
Section~\ref{notation}, we introduce basic notation, the appropriate
function spaces and formulate the problem in the \newtext{strong} and
variational sense. In Section~\ref{assumptions}, we specify our
assumptions on data and modify structure conditions to obtain a
reasonably simple but still realistic model of hygro-thermal
behavior of concrete at high temperatures due to Ba\v{z}ant
and~Thonguthai~\cite{BaTh1978}. An application of the Rothe method
of discretization in time leads to a coupled system of semilinear
steady-state equations, which (together with the appropriate
boundary conditions) form a semilinear elliptic boundary value
problem, formulated in the form of operator equation in appropriate
function spaces. The existence result for this problem in space
$W^{1,p}(\Omega)^2$ with $p \in (2,4)$ is proven in
Section~\ref{existence_discrete} using the general theory of
coercive and pseudomonotone operators in Banach spaces. Next, the
problem is resolved using the finite element method as presented in
Section~\ref{FEM_impl}. In Section \ref{example}, numerical
experiments are performed to investigate the moisture migration,
temperature distribution and pore pressure build up in the model of
concrete specimen exposed to fire, including a simple engineering
approach to study the spalling phenomenon.

\section{The coupled model for wet concrete}\label{model}

\subsection{Conservation laws}\label{Conservation_laws}
The heat and mass transport in concrete is governed by the following
system of conservation laws:

\emph{energy conservation equation:}
\begin{equation}\label{energy temp}
\frac{\partial \mathcal{H}(\theta,w)}{\partial t} = -\nabla \cdot
\bfJ_{\theta}(\theta,w,\nabla\theta,\nabla w) + C_w
\bfJ_w(\theta,w,\nabla\theta,\nabla w)\cdot\nabla\theta;
\end{equation}

\emph{water content conservation equation:}
\begin{equation}\label{conservation_mass}
\frac{\partial \mathcal{M}(\theta,w)}{\partial t} =
-\nabla\cdot\bfJ_w(\theta,w,\nabla\theta,\nabla w).
\end{equation}
The primary unknowns in the balance equations \eqref{energy
temp}--\eqref{conservation_mass} are the temperature $\theta$ and
the water content $w$; $w$ represents the mass of all evaporable
water (free, i.e. not chemically bound) per m$^3$ of concrete.
Further, $\mathcal{H}$ and $\mathcal{M}$ represent the amount of
(internal) energy and the amount of free water, respectively, in $1$
m$^3$ of concrete, $\bfJ_{\theta}$ is the heat flux, $C_w$ the
isobaric heat capacity of bulk (liquid) water and $\bfJ_w$ the
\newtext{humidity} flux.

\subsection{Constitutive relationships for heat and moisture flux}
Following \cite{BaTh1978}, the heat flux $\bfJ_{\theta}$ arises due
to the temperature gradient (Fourier's law) and due to the water
content gradient (Dufour flux)
\begin{equation}\label{heat_flux}
\bfJ_{\theta}(\theta,w,\nabla\theta,\nabla w) =  - D_{\theta
w}(\theta,w) \nabla w - D_{\theta \theta}(\theta,w) \nabla \theta
\end{equation}
and the flux of humidity $\bfJ_w$ consists of the flux due to the
humidity gradient (Fick's law) and due to the temperature gradient
(Soret flux)
\begin{equation}\label{moisture_flux}
\bfJ_{w}(\theta,w,\nabla\theta,\nabla w) = - D_{ww}(\theta,w) \nabla
w - D_{w\theta}(\theta,w) \nabla \theta,
\end{equation}
where $D_{\theta w}$, $D_{\theta \theta}$, $D_{ww}$ and
$D_{w\theta}$ are continuous diffusion coefficient functions
depending non-linearly on $\theta$ and $w$.

\subsection{Boundary and initial conditions}

To complete the introduction of the model, let us specify the boundary and
initial conditions on $\theta$ and $w$. The humidity flux across the boundary is
quantified by the Newton law:
\begin{equation}
\bfJ_{w}(\theta,w,\nabla\theta,\nabla w)\cdot\bfn = \gamma_c (w -
w_{\infty}),
\end{equation}
where the right hand side represents the humidity dissipated into the
surrounding medium with water content $w_\infty$, specified in terms of the film
coefficient $\gamma_c$. As for the heat flux, we shall distinguish the
convective and radiation boundary conditions given by
\begin{eqnarray}
\bfJ_{\theta}(\theta,w,\nabla\theta,\nabla w)\cdot\bfn &=&
\alpha_c(\theta-\theta_{\infty}), \label{neumann a}
\\
\bfJ_{\theta}(\theta,w,\nabla\theta,\nabla w)\cdot\bfn &=&
\alpha_c(\theta-\theta_{\infty})+e\sigma
(\theta|\theta|^3-\theta^4_{\infty}), \label{radiative a}
\end{eqnarray}
respectively, in which $\alpha_c$ designates the film coefficient for
the heat transfer, and $\theta_\infty$ is temperature of the environment.
\newtext{The last expression in Eq.~\eqref{radiative a} expresses the radiative
contribution to the heat flux, quantified by the Stefan-Boltzmann law in terms
of the relative surface emissivity $e$ and the Stefan-Boltzmann constant
$\sigma$ and the temperature difference $( \theta^4 -
\theta_\infty^4)$.\footnote{\newtext{Replacing the term $\theta^4$ with
$\theta|\theta|^3$ is essential later in the proof of Theorem~2.}}} The initial
conditions are set as follows:
\begin{equation}\label{initial_cond_0}
\theta(0) = \theta_0, \qquad w(0) = w_0 .
\end{equation}
Here, $\theta_0$ and $w_0$ represent the initial distributions of the primary
unknowns $\theta$ and $w$, respectively.
\section{Notation and formulation of the problem}
\label{notation}

Vectors, vector functions and operators acting on vector functions
are denoted by~boldface letters. Throughout the paper, we will
always use positive constants $c$, $c_1$, $c_2$, $\dots$, which are
not specified and which may differ from line to line. For an
arbitrary $r\in [1,+\infty]$, $L^r(\Omega)$ denotes the usual
Lebesgue space equipped with the norm $\|\cdot\|_{L^r(\Omega)}$, and
$W^{k,p}(\Omega)$, $k\geq 0$, $p\in [1,+\infty]$, denotes the usual
Sobolev space with the norm $\|\cdot\|_{W^{k,p}(\Omega)}$. Let $X$
be a Banach space. By $C([0,T],X)$ we \newtext{denote} the space of all
continuous functions $\varphi:[0,T]\rightarrow X$. Throughout the
paper $p'=p/(p-1)$, \newtext{$p > 1$}, denotes the conjugate exponent to $p$.
$\bfphi'(t)$ indicates the partial derivative with respect to time;
we also write $\bfphi'(t) = \bfphi_t$.

We consider a mixed initial--boundary value problem for a general
model of the coupled heat and mass flow in a two-dimensional domain
$\Omega\subset \mathbb{R}^2$ with a Lipschitz boundary
$\partial\Omega$, which consists of non-intersecting pieces
$\Gamma_R$ and $\Gamma_N$, $\partial \Omega = \overline{\Gamma_R}
\cup \overline{\Gamma_N}$. $\Gamma_R$ re  presents the part of the
boundary which is exposed to fire, whereas the other part denoted by
$\Gamma_N$ is exposed to atmosphere. Let $T>0$ be the fixed value of
the time horizon, $Q_{T}=\Omega\times(0,T)$,
$\Gamma_{RT}=\Gamma_R\times(0,T)$ and
$\Gamma_{NT}=\Gamma_N\times(0,T)$. The \newtext{strong} formulation of our
problem is as follows:
\begin{align}
\mathcal{H}_t &= -\nabla \cdot \bfJ_{\theta} + C_w \bfJ_w
\cdot\nabla \theta   &&{\rm in} \; Q_{T},\label{eq1}
\\
\mathcal{M}_t&= -\nabla \cdot \bfJ_w  && {\rm in}\;   Q_{T},
\\
\bfJ_{w}\cdot\bfn &= \gamma_c(w-w_{\infty})  &&{\rm on} \;
\Gamma_{NT} \cup \Gamma_{RT},
\\
\bfJ_{\theta} \cdot \bfn  &= \alpha_c(\theta-\theta_{\infty}) &&
{\rm on} \;  \Gamma_{NT},
\\
\bfJ_{\theta} \cdot \bfn &= \alpha_c(\theta-\theta_{\infty})+
e\sigma (|\theta|^3\theta-\theta^4_{\infty}) &&{\rm on} \;
\Gamma_{RT},
\\
\theta(0)&= \theta_0 &&{\rm in} \;  \Omega,
\\
w(0)&= w_0 &&{\rm in} \;  \Omega.\label{eq10}
\end{align}
Here we assume that all functions are smooth enough. Now we can
formulate the problem in the variational sense. Suppose that
$[\theta_{\infty}(t),w_{\infty}(t)]\in C(0,T)^2$  and $[\theta_0,
w_0]\in W^{1,r}(\Omega)^2$, $r>2$. Find a pair $[\theta,w]\in
C([0,T];W^{1,r}(\Omega)^2)$ such that
\begin{multline*}
-\int_0^T \left\langle \mathcal{H},\phi'_1(t) \right\rangle +
\left\langle \mathcal{M},\phi'_2(t)\right\rangle {\rm d}t
-\int_{Q_T}\bfJ_{\theta}\cdot\nabla\phi_1 +
\bfJ_{w}\cdot\nabla\phi_2{\rm d}Q_T
\\
-\int_{Q_T}C_w \bfJ_w \cdot\nabla \theta\phi_1 \;{\rm d}Q_T
+\int_{\Gamma_{RT}} e\sigma (|\theta|^3\theta-\theta^4_{\infty})
\phi_1 \;{\rm d}S_T
\\
+\int_{\Gamma_{RT}\cup\Gamma_{NT}} \alpha_c(\theta-\theta_{\infty})
\phi_1 \;{\rm d}S_T
+\int_{\Gamma_{RT}\cup\Gamma_{NT}} \gamma_c(w-w_{\infty})  \phi_2
\;{\rm d}S_T = 0
\end{multline*}
holds for all test functions $\bfphi=[\phi_1,\phi_2]\in
C_0^{\infty}(  [0,T];C^{\infty}(\overline{\Omega})^2)$ and
\begin{equation}\label{initial_cond}
\theta(0)=\theta_0 \;\textmd{ and }\; w(0)=w_0  \quad \textmd{ in }
\Omega.
\end{equation}
Such pair $[\theta,w]$ is called the variational solution to the
system \eqref{eq1}--\eqref{eq10}.

\begin{rmk}
\newtext{To the best of our knowledge, there are no existence results
for the presented model available.}
\end{rmk}

\section{Existence of the approximate solution to the Ba\v{z}ant's model}

\subsection{Structural conditions and assumptions on physical parameters}
\label{assumptions}

\begin{itemize}
\item[$A_1$]
In \cite{BaTh1978}, Ba\v{z}ant and Thonguthai expressed the time
variation of amounts $\mathcal{H}$ and $\mathcal{M}$ as follows:
\begin{eqnarray}\label{amount_H}
\frac{\partial \mathcal{H}}{\partial t} &:=& \rho_s C_s
\frac{\partial \theta}{\partial t}-h_{d}\frac{\partial w_d}{\partial
t}-h_{\alpha}\frac{\partial w}{\partial t},
\\
\frac{\partial \mathcal{M}}{\partial t} &:=& \frac{\partial
w}{\partial t}-\frac{\partial w_d}{\partial t}.
\end{eqnarray}
Here $\rho_s$ and $C_s$, respectively, are the mass density and the
isobaric heat capacity of solid microstructure (excluding hydrate
water), $w_d$ represents the total mass of the free water released
in the pores by drying. $h_{\alpha}$ denotes the enthalpy of
evaporation per unit mass and $h_{d}$ denotes the enthalpy of
dehydration per unit mass.
\item[$A_2$] We assume that the parameters $\rho_s$, $C_w$, $h_d$, $\alpha_c$, $\beta_c$,
$\sigma$ and $e$ are real positive constants.
\item[$A_3$] The functions $C_s=C_s(\theta)$ and $h_{\alpha}=h_{\alpha}(\theta)$ are positive
continuous functions, $w_d=w_d(\theta)$ is positive increasing
function belonging to $W^{1,\infty}(\mathbb{R})$,
$\theta_{\infty}(t)$ and $w_{\infty}(t)$ are given continuous
functions of time and $\theta_0, w_0 \in W^{1,r}(\Omega)$, $r>2$.

\item[$A_4$] Water content $w$ is connected with temperature $T$ and pore
pressure $P$ via a so-called sorption isotherm $w=\Phi(\theta,P)$, which has to
be determined experimentally for each type of concrete. We assume $\Phi$ to be a
continuous function such that $\Phi(\xi_1,\xi_2)\geq 0$ for $\bfxi \in
\mathbb{R}_+^2$ and $\Phi=0$ otherwise.
\item[$A_5$] Following \cite{BaTh1978}, we consider the cross effects
to be negligible. This leads to simple phenomenological relations introduced by
Ba\v{z}ant and Thonguthai in the form
\begin{equation}
\label{heat_flux_bazant} \bfJ_{\theta}:=-\lambda_c(\theta) \nabla
\theta \quad \textmd{ and  } \quad  \bfJ_{w}:=
-\frac{\kappa(\theta,P)}{g} \nabla P,
\end{equation}
where the thermal conductivity $\lambda_c$ and permeability $\kappa$ are assumed
to be positive continuous functions of their arguments and $g$ is the
gravitational acceleration.

\end{itemize}

\subsection{Solutions to the discretized problem}\label{existence_discrete}
Incorporating the relations
\eqref{amount_H}--\eqref{heat_flux_bazant} into the system
\eqref{eq1}--\eqref{eq10} we get the modified Ba\v{z}ant-Thonguthai
model with primary unknowns $w$, $\theta$ and $P$ consisting of

\noindent \emph{conservation laws:}
\begin{align}
\frac{\partial w}{\partial t}  &= \nabla
\cdot\left(\frac{\kappa(\theta,P)}{g} \nabla P \right) +
\frac{\partial w_d(\theta)}{\partial t} &&{\rm in} \;
Q_{T},\label{eq100}
\\
\rho_s C_s(\theta)\frac{\partial \theta}{\partial
t}-h_{\alpha}(\theta)\frac{\partial w}{\partial t} &= \nabla \cdot
(\lambda_c(\theta,P) \nabla \theta)
\nonumber\\
& \quad - C_w \frac{\kappa(\theta,P)}{g}\nabla P \cdot\nabla \theta
+ h_{d}\frac{\partial w_d(\theta)}{\partial t} && {\rm in}\; Q_{T};
\end{align}
\emph{state equation of pore water:}
\begin{align}
w - \Phi(P,\theta) = 0 \quad \textmd{  in } Q_{T}; \label{eq100c}
\end{align}
\emph{radiation boundary conditions:}
\begin{align}
- \lambda_c(\theta,P) \nabla \theta \cdot \bfn  &=
\alpha_c(\theta-\theta_{\infty})+ e\sigma
(|\theta|^3\theta-\theta^4_{\infty}) &&{\rm on} \; \Gamma_{RT};
\end{align}
\emph{Neumann boundary conditions:}
\begin{align}
- \lambda_c(\theta,P) \nabla \theta \cdot \bfn  &=
\alpha_c(\theta-\theta_{\infty}) && {\rm on} \;  \Gamma_{NT},
\\
-\frac{\kappa(\theta,P)}{g}\nabla P\cdot\bfn &=
\beta_c(P-P_{\infty}) &&{\rm on} \; \Gamma_{NT} \cup \Gamma_{RT};
\end{align}
and \emph{initial conditions:}
\begin{equation}\label{eq102}
P(0) = P_0  \quad \textmd{ and } \quad  \theta(0) = \theta_0 \qquad
\textmd{ in } \Omega.
\end{equation}

Let $0=t_0<t_1<\dots<t_N=T$ be an equidistant partitioning of time interval
$[0;T]$ with step $\Delta t$. Set a fixed integer $n$  such that $0\leq n <N$.
In what follows we abbreviate $\varphi(\bfx,t_n)$ by $\varphi_n$ for any
function $\varphi$. The time discretization of the continuous model is
accomplished through a semi-implicit difference scheme
\begin{equation}\label{eq:semi-implicit-1}
\frac{w_{n+1}-w_n}{\Delta t} =\nabla
\cdot\left(\frac{\kappa(\theta_n,P_n)}{g} \nabla P_{n+1}\right) +
\frac{w_d(\theta_{n+1})-w_d(\theta_n)}{\Delta t},
\end{equation}
\begin{multline}\label{eq:semi-implicit-2}
\rho_s C_s(\theta_n)\frac{\theta_{n+1}-\theta_n}{\Delta t}-
h_{\alpha}(\theta_n)\frac{w_{n+1}-w_n}{\Delta t} = \nabla \cdot
(\lambda_c(\theta_n,P_n) \nabla \theta_{n+1})
\\
- C_w \frac{\kappa(\theta_n,P_n)}{g}\nabla P_n \cdot\nabla \theta_n
+ h_{d}\frac{w_d(\theta_{n+1})-w_d(\theta_n)}{\Delta t}.
\end{multline}
Here, we assume that the functions $\theta_n$, $w_n$ and $P_n$ are
known. In what follows we study the problem of existence of the
solution $\theta_{n+1}$, $w_{n+1}$ and $P_{n+1}$. Incorporating the
relation \eqref{eq100c} into the system \eqref{eq100}--\eqref{eq102}
we can eliminate the unknown field $w$ and consider the problem with
only two unknowns $\theta$ and $P$. Consequently, the existence of
$w_{n+1}$ follows from the existence of $\theta_{n+1}$ and $P_{n+1}$
by the relation \eqref{eq100c}. For the sake of simplicity we denote
$[\pi,\tau]:=[P_{n+1},\theta_{n+1}]$. Let us put
$\widetilde{\kappa}(\bfx) ={\kappa(\theta_n(\bfx),P_n(\bfx))}/{g}$,
$\widetilde{\lambda}_c(\bfx)= \lambda_c(\theta_n(\bfx),P_n(\bfx))$,
$\Phi_n=\Phi(P_n,\theta_n)$ and introduce the functions
\begin{eqnarray}
R_1(\bfx,\pi,\tau)&=& \frac{1}{\Delta t}\Phi(\pi,\tau) -
\frac{1}{\Delta t} w_d(\tau),
\label{def_R1}
\\
R_2(\bfx,\pi,\tau)&=& \frac{1}{\Delta t}\rho_s C_s(\theta_n)\tau-
\frac{1}{\Delta t}h_{d} w_d(\tau) -\frac{1}{\Delta
t}h_{\alpha}(\theta_n)\Phi(\pi,\tau), \label{def_R2}
\\
F_1(\bfx)&=&  \frac{1}{\Delta t}\Phi_n - \frac{1}{\Delta t}
w_d(\theta_n),\label{def_F1}
\\
F_2(\bfx)&=& \frac{1}{\Delta t}\rho_s
C_s(\theta_n)\theta_n-\frac{h_{d}}{\Delta t}
w_d(\theta_n)-\frac{h_{\alpha}(\theta_n)}{\Delta
t}\Phi_n-C_w\widetilde{\kappa}(\bfx)\nabla P_n \cdot\nabla \theta_n.
\nonumber
\\\label{def_F2}
\end{eqnarray}
Obviously, we have to solve, successively for $n=0,\ldots,N-1$, the
following semilinear system with primary unknowns $[\pi,\tau]$
\begin{equation}\label{problem_eq1}
\left\{
\begin{array}{rclll}
- \nabla \cdot({\widetilde{\kappa}(\bfx)} \nabla \pi)
+R_1(\bfx,\pi,\tau) &=& F_1(\bfx)
 & {\rm in } & \Omega,\\
-\nabla \cdot (\widetilde{\lambda}_c(\bfx) \nabla
\tau)+R_2(\bfx,\pi,\tau)
 &=& F_2(\bfx) &  {\rm in } & \Omega
\end{array}\right.
\end{equation}
and with the boundary conditions
\begin{equation}\label{bc}
\left\{
\begin{array}{rclll}
- \widetilde{\lambda}_c(\bfx) \nabla \tau\cdot\bfn&=& \alpha_c(\tau
-\theta_{\infty,n})+ e\sigma (|\tau|^3\tau-\theta^4_{\infty,n})&
{\rm on }&\Gamma_{R},
\\
-\widetilde{\kappa}(\bfx) \nabla \pi \cdot \bfn &=& \beta_c(\pi
-P_{\infty,n})&{\rm on }&\Gamma_{N} \cup \Gamma_{R},
\\
- \widetilde{\lambda}_c(\bfx) \nabla \tau \cdot \bfn
&=&\alpha_c(\tau -\theta_{\infty,n})& {\rm on} & \Gamma_{N}.
\end{array}\right.
\end{equation}
\begin{dfn}
The pair $[\pi,\tau]\in W^{1,r}(\Omega)^2$, $r\geq 2$, is called a
variational solution to the system \eqref{problem_eq1}--\eqref{bc}
iff
\begin{multline}\label{weak_form}
\int_{\Omega} {\widetilde{\kappa}(\bfx)}\nabla\pi\cdot\nabla v_{\pi}
+\widetilde{\lambda}_c(\bfx)\nabla\tau\cdot\nabla v_{\tau}\,{\rm
d}\bfx + \int_{\Omega}
R_1(\bfx,\pi,\tau)\,v_{\pi}+R_2(\bfx,\pi,\tau)\,v_{\tau}\,{\rm
d}\bfx
\\
+ \int_{\Gamma_R\cup\Gamma_N} \!\!\! \beta_c\pi\,v_{\pi}\,{\rm
d}\bfS + \int_{\Gamma_R\cup\Gamma_N} \!\!\!\alpha_c\tau v_{\tau}{\rm
d}\bfS + \!\!\int_{\Gamma_R}\!\!\!e\sigma |\tau|^3\tau v_{\tau} {\rm
d}\bfS
\\
= \int_{\Omega}\! F_1(\bfx) v_{\pi}+ F_2(\bfx) v_{\tau}{\rm d}\bfx +
\int_{\Gamma_R\cup\Gamma_N}\!\!\! \beta_c P_{\infty,n}v_{\pi}+
\alpha_c \theta_{\infty,n}v_{\tau}{\rm d}\bfS+ \int_{\Gamma_R}\!\!\!
e\sigma \theta^4_{\infty,n} v_{\tau}{\rm d}\bfS
\end{multline}
holds for every $[v_{\tau},v_{\pi}]\in W^{1,r'}(\Omega)^2$,
$r'=r/(r-1)$.
\end{dfn}
The main result of this section is the following
Theorem~\ref{existence_thm} and Corollary~\ref{main_cor}.
\begin{thm}\label{existence_thm}
Assume that $[P_n,\theta_n]\in W^{1,p}(\Omega)^2$ with some fixed
$p\in(2,4)$ is known and let $A_2$--$A_5$ be satisfied. Then there
exists the variational solution $[\pi,\tau]\in W^{1,p}(\Omega)^2$ to
the system \eqref{problem_eq1}--\eqref{bc}.
\end{thm}
\begin{pf}
In order to prove Theorem~\ref{existence_thm}, it is convenient to
define the operator $\mathcal{T}:W^{1,p}(\Omega)^2\rightarrow
W^{-1,p}(\Omega)^2$ by
\begin{eqnarray}\label{def_T}
\langle
\mathcal{T}([\pi,\tau]),[v_{\pi},v_{\tau}]\rangle&=&\int_{\Omega}
({\widetilde{\kappa}(\bfx)} \nabla \pi) \cdot \nabla v_{\pi} \,{\rm
d}\bfx +\int_{\Omega} R_1(\bfx,\pi,\tau) \; v_{\pi} \,{\rm d}\bfx
\nonumber \\
&& + \int_{\Omega} (\widetilde{\lambda}_c(\bfx) \nabla \tau)  \cdot
\nabla v_{\tau} \,{\rm d}\bfx +\int_{\Omega} R_2(\bfx,\pi,\tau) \;
v_{\tau} \,{\rm d}\bfx
\nonumber\\
&& + \! \int_{\Gamma_R\cup\Gamma_N}\!\!\!\!\! \beta_c \pi v_{\pi}
{\rm d}\bfS + \!\int_{\Gamma_R\cup\Gamma_N} \!\!\!\!\!\!\alpha_c\tau
v_{\tau} {\rm d}\bfS + \!\!\int_{\Gamma_R}\!\!\!e\sigma |\tau|^3\tau
v_{\tau} \,{\rm d}\bfS
\nonumber\\
\end{eqnarray}
and the functional $\bff \in W^{-1,p}(\Omega)^2$ by
\begin{eqnarray}\label{def_f}
\langle\bff,[v_{\pi},v_{\tau}]\rangle&=& \int_{\Omega}\! F_1(\bfx)
v_{\pi}{\rm d}\bfx + \int_{\Omega}\! F_2(\bfx) v_{\tau}{\rm d}\bfx +
\int_{\Gamma_R\cup\Gamma_N}\!\!\! \beta_c P_{\infty,n}v_{\pi}{\rm
d}\bfS
\nonumber \\
&&  + \int_{\Gamma_R\cup\Gamma_N}\!\!\!  \alpha_c \theta_{\infty,n}
v_{\tau} {\rm d}\bfS  +  \int_{\Gamma_R} e\sigma \theta^4_{\infty,n}
v_{\tau} {\rm d}\bfS
\end{eqnarray}
for all $[v_{\pi},v_{\tau}]\in W^{1,p'}(\Omega)^2$. Since we assume
$[P_n,\theta_n]\in W^{1,p}(\Omega)^2$ with some fixed $p\in(2,4)$,
we have the following estimate for the convective term
$C_w\widetilde{\kappa}(\bfx)\nabla P_n \cdot\nabla \theta_n$
\begin{equation}
\int_{\Omega}\! \left(C_w\widetilde{\kappa}(\bfx)\nabla P_n
\cdot\nabla \theta_n\right) v_{\tau}{\rm d}\bfx \leq c_1
\|P_n\|_{W^{1,p}(\Omega)}\|\theta_n\|_{W^{1,p}(\Omega)}
\|v_{\tau}\|_{W^{1,p'}(\Omega)}
\end{equation}
for all $v_{\tau} \in W^{1,p'}(\Omega)$. One can prove in the
similar way that the other integrals in \eqref{def_f} are finite and
the functional $\bff$ is well-defined. The variational problem can
now be treated as a single operator equation
$\mathcal{T}([\pi,\tau])=\bff$. Obviously, $\bff\in
W^{-1,p}(\Omega)^2$ implies $\bff\in W^{-1,2}(\Omega)^2$. \newtext{
First of all we prove that for a given $\bfh\in W^{-1,2}(\Omega)^2$
there exists $[\pi,\tau]\in W^{1,2}(\Omega)^2$: the solution of the
equation $\mathcal{L}([\pi,\tau])=\bfh$ with
$\mathcal{L}:W^{1,2}(\Omega)^2\rightarrow W^{-1,2}(\Omega)^2$
defined by \eqref{def_T} (substituting $\mathcal{L}$ instead of
$\mathcal{T}$) for all $[v_{\pi},v_{\tau}]\in W^{1,2}(\Omega)^2$. }
\begin{lem}\label{T_bounded}
$\newtext{\mathcal{L}}:W^{1,2}(\Omega)^2\rightarrow W^{-1,2}(\Omega)^2$
is bounded.
\end{lem}
\begin{pf}
Test \eqref{def_T} by $[\pi,\tau]\in W^{1,2}(\Omega)^2$. Take into
account $A_2$--$A_4$ to get
\begin{eqnarray*}
\langle \newtext{\mathcal{L}}([\pi,\tau]),[\pi,\tau]\rangle &\leq&
c_1\|\pi\|_{W^{1,2}(\Omega)}^2+c_2\|\tau\|_{W^{1,2}(\Omega)}^2
\\
&& + c_3\|\pi\|_{L^{2}(\Omega)}^2+c_4\|\tau\|_{L^{2}(\Omega)}^2
\\
&& + \beta_c\|\pi\|_{L^{2}(\partial\Omega)}^2+
\alpha_c\|\tau\|_{L^{2}(\partial\Omega)}^2 +
e\sigma\|\tau\|_{L^{5}(\partial\Omega)}^5
\\
&\leq&  c_5 \|[\pi,\tau]\|_{W^{1,2}(\Omega)^2}^2 +
e\sigma\|\tau\|_{L^{5}(\partial\Omega)}^5 .
\end{eqnarray*}
Due to the trace theorem \cite{Necas1983} there exists a constant
$c_{tr}$ such that
\begin{equation*}
\|v\|_{L^q(\partial\Omega)} \leq c_{tr}\|v\|_{W^{1,2}(\om)} \textmd{
for all } v\in W^{1,2}(\om), \; q \geq 1.
\end{equation*}
Hence $\newtext{\mathcal{L}}$ is bounded. $\square$
\end{pf}
\begin{lem}\label{T_coercive}
$\newtext{\mathcal{L}}:W^{1,2}(\Omega)^2\rightarrow W^{-1,2}(\Omega)^2$
is coercive.
\end{lem}
\begin{pf}
$A_3$, $A_4$  and the Young inequality yield
\begin{eqnarray}\label{R1_coercivity}
R_1(\bfx,\xi_1,\xi_2)\xi_1 &=& \left(\frac{1}{\Delta
t}\Phi(\xi_1,\xi_2) - \frac{1}{\Delta t}w_d(\xi_2)\right)\xi_1
\nonumber\\
&=& \frac{1}{\Delta t}\Phi(\xi_1,\xi_2)\xi_1 - \frac{1}{\Delta
t}w_d(\xi_2)\xi_1
\nonumber\\
& \geq&-\eta { \xi_1^2 } - c(\eta)\left(\frac{1}{\Delta
t}w_d(\xi_2)\right)^2
\end{eqnarray}
for every $\xi_1,\xi_2\in \mathbb{R}$ and arbitrary $\eta>0$.
Further, $A_3$, $A_4$  and the Young inequality  yield the existence
of a positive function $g_1$ and a non-negative function $g_2$ (both
of spatial variable $\bfx$) such that
\begin{eqnarray}\label{R2_coercivity}
R_2(\bfx,\xi_1,\xi_2)\xi_2 &=& \frac{1}{\Delta t}\rho_s
C_s(\theta_n)\xi^2_2 -\frac{1}{\Delta t}\left(h_{d}  w_d(\xi_2)
 +h_{\alpha}(\theta_n)\Phi(\xi_1,\xi_2) \right) \xi_2
\nonumber\\
 & \leq & g_1(\bfx)\xi^2_2 - g_2(\bfx) \qquad \forall \bfx \in
 \Omega,\; \forall \, [\xi_1,\xi_2]\in \mathbb{R}^2.
\end{eqnarray}
Now \eqref{def_T}, \eqref{R1_coercivity}, \eqref{R2_coercivity}, the
embedding $W^{1,2}(\Omega)\hookrightarrow L^2(\Omega)$ and the
Friedrichs inequality imply
\begin{eqnarray}\label{T_coercivity}
\langle \newtext{\mathcal{L}}([\pi,\tau]),[\pi,\tau]\rangle &=&
\int_{\Omega} {\widetilde{\kappa}(\bfx)} |\nabla \pi|^2  \,{\rm
d}\bfx + \int_{\Gamma_R\cup\Gamma_N}  \!\!\!  \beta_c |\pi|^2 \,{\rm
d}\bfS
\nonumber \\
&& + \int_{\Omega} \widetilde{\lambda}_c(\bfx) |\nabla \tau|^2
\,{\rm d}\bfx +  \!\! \int_{\Gamma_R\cup\Gamma_N} \!\!\!\!\!\!
\alpha_c|\tau|^2 \,{\rm d}\bfS +  \!\! \int_{\Gamma_R} \!\!\!
e\sigma|\tau|^3\tau^2\,{\rm d}\bfS
\nonumber\\
&&  + \int_{\Omega} R_1(\bfx,\pi,\tau)\;\pi \,{\rm d}\bfx
+\int_{\Omega} R_2(\bfx,\pi,\tau)\;\tau\,{\rm d}\bfx
\nonumber \\
&\geq& \!
c_1\|\pi\|_{W^{1,2}(\Omega)}^2\!+\!c_2\|\tau\|_{W^{1,2}(\Omega)}^2
\!-\! \eta \|\pi\|_{L^{2}(\Omega)}^2 \!+\!
c_3\|\tau\|_{L^{2}(\Omega)}^2 \!-\! c_4
\nonumber \\
&\geq&  c_5 \|[\pi,\tau]\|_{W^{1,2}(\Omega)^2}^2 - c_6
\end{eqnarray}
with some positive constants $c_1,\dots,c_6$ and choosing $\eta$
sufficiently small. $\square$
\end{pf}
\begin{lem}\label{T_pseudomonotone}
$\newtext{\mathcal{L}}:W^{1,2}(\Omega)^2\rightarrow W^{-1,2}(\Omega)^2$
is pseudomonotone.
\end{lem}
\begin{pf}
Obviously, since $\widetilde{\kappa}>0$ and
$\widetilde{\lambda}_c>0$ in $\Omega$, the inequality
\begin{equation*}
{\widetilde{\kappa}(\bfx)}(\xi_1-\xi'_1)^2
+\widetilde{\lambda}_c(\bfx)(\xi_2-\xi'_2)^2
>0
\end{equation*}
holds for all $\bfx \in \Omega$ and for all $[\xi_1,\xi_2],[\xi'_1,\xi'_2]\in
\mathbb{R}^2$, $[\xi_1,\xi_2]\neq[\xi'_1,\xi'_2]$. $\square$
\end{pf}
\begin{crl}\label{existence_p2}
The smoothness assumptions on $\widetilde{\kappa}$,
$\widetilde{\lambda}_c$, $R_1$ and $R_2$ (the smoothness of $R_1$
and $R_2$ follows from $A_2$--$A_4$ and \eqref{def_R1} and
\eqref{def_R2}) and Lemma \ref{T_bounded}--\ref{T_pseudomonotone}
imply that $\newtext{\mathcal{L}}$ is continuous, bounded, coercive and
pseudomonotone. Now \cite[Theorem 3.3.42]{Necas1983} yields the
existence of the solution $[\pi,\tau]\in W^{1,2}(\Omega)^2$ to the
equation $\newtext{\mathcal{L}}([\pi,\tau])=\newtext{\bfh}$ for every
$\newtext{\bfh} \in W^{-1,2}(\Omega)^2$.
\end{crl}
To get higher regularity results go back and consider $\bff\in
W^{-1,p}(\Omega)^2 \subset W^{-1,2}(\Omega)^2$ with some $p\in(2,4)$
and rewrite the system \eqref{problem_eq1}--\eqref{bc} in the form
\begin{equation}\label{eq50}
\left\{
\begin{array}{rclll}
- \nabla \cdot({\widetilde{\kappa}(\bfx)} \nabla \pi) &=&
F_1(\bfx)-R_1(\bfx,\pi,\tau) &\mathrm{in}& \Omega,
\\
-\nabla \cdot (\widetilde{\lambda}_c(\bfx) \nabla \tau) &=&
F_2(\bfx)-R_2(\bfx,\pi,\tau)&{\rm in}& \Omega,
\\
- \widetilde{\lambda}_c(\bfx) \nabla \tau\cdot\bfn&=& \alpha_c(\tau
-\theta_{\infty,n})+ e\sigma (|\tau|^3\tau-\theta^4_{\infty,n})
&\mathrm{on}&\Gamma_{R},
\\
-\widetilde{\kappa}(\bfx) \nabla \pi \cdot \bfn &=& \beta_c(\pi
-P_{\infty,n})&{\rm on}&\Gamma_{R} \cup \Gamma_{N},
\\
- \widetilde{\lambda}_c(\bfx) \nabla \tau \cdot \bfn
&=&\alpha_c(\tau -\theta_{\infty,n})&{\rm on}& \Gamma_{N}.
\end{array}\right.
\end{equation}
It is easy to verify that for the functional $\bff\in
W^{-1,p}(\Omega)^2 \subset W^{-1,2}(\Omega)^2$ defined by
\eqref{def_f} and for the weak solution $[\pi,\tau]\in
W^{1,2}(\Omega)^2$ (whose existence is ensured by Corollary
\ref{existence_p2}), the functional $\bfG\in W^{-1,p}(\Omega)^2$
given by
\begin{multline}\label{functional_G}
\langle\bfG,[v_{\pi},v_{\tau}]\rangle=\int_{\Omega}\!
\left(F_1-R_1\right) v_{\pi} + \left(F_2-R_2\right) v_{\tau}{\rm
d}\bfx + \int_{\Gamma_R\cup\Gamma_N}\!\!\! \alpha_c(\tau
-\theta_{\infty,n}) v_{\tau}{\rm d}\bfS
\\
+ \int_{\Gamma_R\cup\Gamma_N}\!\!\! \beta_c(\pi
-P_{\infty,n})v_{\pi}{\rm d}\bfS + \int_{\Gamma_R} \left(e\sigma
(|\tau|^3\tau-\theta^4_{\infty,n})\right) v_{\tau}{\rm d}\bfS
\end{multline}
for every $[v_{\pi},v_{\tau}]\in W^{1,p'}(\Omega)^2$ is well
defined. It is known (see \cite{KufSan1987,MazRoss2003}) that for
given $\bfG\in W^{-1,p}(\Omega)^2$ defined by \eqref{functional_G}
with $p\in(2,4)$ the Neumann problem for the elliptic system
\eqref{eq50} (where the right hand side is represented by $\bfG$)
possess the solution $[\pi,\tau]\in W^{1,p}(\Omega)^2$. This
completes the proof of Theorem~\ref{existence_thm}.
\end{pf}

\begin{crl}\label{main_cor}
Following Theorem~\ref{existence_thm}, $[P_{k},\theta_{k}]\in
W^{1,p}(\Omega)^2$  yields $[P_{k+1},\theta_{k+1}]\in
W^{1,p}(\Omega)^2$ with any $p\in(2,4)$. Since we suppose
$[P_{0},\theta_{0}]\in W^{1,r}(\Omega)^2$ with $r > 2$, we can
conclude, that $[P_{n},\theta_{n}]\in W^{1,p}(\Omega)^2$
successively for $n=0,\ldots,N-1$ for any $p\in(2,r)$ if $r<4$ and
$p\in(2,4)$ if $r\geq 4$. \newtext{Note that this solution needs not to be
unique.}
\end{crl}

\section{Numerical results}

\subsection{Finite element implementation}\label{FEM_impl}

Consider a polygonal approximation $\Omega^h$ to $\Omega$, defined
by an admissible quadrilateral partition $\mathcal{\newtext{Q}}^h = \{
\newtext{Q}_1, \newtext{Q}_2, \ldots, \newtext{Q}_{N_e} \}$ such that
$\overline{\Omega}^h = \cup_{e=1}^{N_e} \overline{\newtext{Q}}_e$, $\Omega^h
\subseteq \Omega$ and every element $\newtext{Q}_e$ has diameter at most $2h$.
$N_n$ is used to denote the number of nodes of the mesh and $\Gamma_R^h$
or $\Gamma_N^h$ stand for the part of the approximate boundary
$\Gamma^h = \partial \Omega^h$ where the radiation and convection
boundary conditions are prescribed. We associate with
$\mathcal{\newtext{Q}}^h$ a finite-dimensional space \newtext{of piecewise
bi-linear basis functions}~(recall that $p \in (2, 4)$)
\begin{eqnarray}
S^h
& = &
\Bigl\{
v \in C^0\left( \overline{\Omega}^h \right)
:
v|_{\newtext{Q}_e} \in \mathcal{P}_2
\text{ and restriction to each edge of } \partial
\newtext{Q}_e
\nonumber \\
&&
\text{ belongs to } \mathcal{P}_1
\text{ for } e = 1, 2, \ldots, N_e
\Bigr\} \subset
W^{1,p}( \Omega^h ),
\end{eqnarray}
where $\mathcal P_s$ denotes the set of polynomials of degree $\leq
s$, cf.~\cite[Section~5]{Braess2007}.

From the implementation point of view, it is more convenient to
derive the numerical scheme by considering all three unknowns $(w,
\theta, P)$ instead of the reduced version~\eqref{weak_form};
the approximate solution $(w^h_{n+1}, \theta^h_{n+1}, P^h_{n+1} )
\in \big[S^h\big]^3$ is thus provided by the weak form of balance
equations~\eqref{eq:semi-implicit-1} and~\eqref{eq:semi-implicit-2}, tested by
$v_w \in S^h$ and $v_P \in S^h$, constrained by the isotherm relation
$A_4$ enforced at the nodes. This leads to a system of non-linear
algebraic equations
\begin{equation}\label{eq:discrete_system}
\frac{1}{\Delta t} \M{C}_n \left( \M{X}_{n+1} - \M{X}_n \right) +
\M{K}_n \M{X}_{n+1} + \M{R}\big( \M{X}_{n+1} \big) = \M{F}_{n+1},
\end{equation}
where e.g. $\M{X}_{n+1} = \left( \M{w}_{n+1}, \M{\theta}_{n+1},
\M{P}_{n+1} \right) \in \mathbb{R}^{3N_n \times 1}$ stores the
unknown nodal values of water content, temperature and pore pressure
at time $t_{n+1}$, respectively. The constant matrices
in~\eqref{eq:discrete_system} exhibit a block structure
\begin{equation}\label{eq:block_structure_linear}
\M{C}_n = \left(
\begin{array}{ccc}
\M{C}^{ww}        & \M{0}                   & \M{0} \\
-\M{C}^{\theta w}_\newtext{n} & \M{C}^{\theta\theta}_n  & \M{0} \\
\M{0}             & \M{0}                   & \M{0}
\end{array}
\right),
\M{K}_n = \left(
\begin{array}{ccc}
\M{0} & \M{0} & \M{K}^{wP}_{n} \\
\M{0} & \M{K}^{\theta\theta}_n & \M{0} \\
\M{0} & \M{0} & \M{0}
\end{array}
\right),
\M{F}_{n+1} = \left(
 \begin{array}{c}
- \M{F}_{n+1}^w \\
 \M{F}_{n+1}^\theta \\
 \M{0}
 \end{array}
\right)
\end{equation}
and the non-linear term reads as
\begin{equation}\label{eq:block_structure_non_linear}
\M{R}\big( \M{X}_{n+1} \big) = \left(
\begin{array}{c}
-\M{R}^{w} ( \M{\theta}_{n+1} ) \\
\M{R}^{\theta} ( \M{\theta}_{n+1} ) \\
\M{w}_{n+1} - \Phi( \M{\theta}_{n+1}, \M{P}_{n+1} )
\end{array}
\right).
\end{equation}

The sub-matrices in~\eqref{eq:block_structure_linear}
and~\eqref{eq:block_structure_non_linear} are obtained by the
assembly of element contributions, e.g.
\begin{eqnarray}
\M{C}_n^{\newtext{\theta w}} = \assembly_{e=1}^{N_e} \M{C}_{n,e}^{\newtext{
\theta w}}, & \displaystyle \M{F}_{n+1}^{w} = \assembly_{e=1}^{N_e}
\M{F}_{n+1,e}^{w}, & \M{R}^{\theta} ( \M{\theta}_{n+1} ) =
\assembly_{e=1}^{N_e} \M{R}_e^\theta \big( \M{\theta}_{n+1,e} \big).
\end{eqnarray}
Here, $\assembly$ is the assembly
operator~\cite[Section~2.8]{Braess2007} and $\theta_{n,e} \in
\mathbb R^{3 \times 1}$ stores the temperature values at the nodes
of the $e$-th element at time $t_n$, related to the local
approximation
\begin{equation}
\theta^h_{n,e}(\bfx )
:=
\theta^h_{n} (\bfx )  \big|_{T_e}
=
\M{N}_e( \bfx ) \M{\theta}_{n,e},
\end{equation}
in which $\M{N}_e : T_e \rightarrow \mathbb R^{1\times 3}$ denotes
the operator of linear basis functions. Analogous relations hold for
the remaining fields.

The individual symmetric positive-definite matrices
$\M{C}_{\bullet}^{\bullet} \in \mathbb R^{3\times 3}$ can now be
expressed as
\begin{subequations}
\begin{eqnarray}
\M{C}^{ww}_e & = & \int_{T_e} \M{N}_e( \bfx )\trn \M{N}_e( \bfx )
\de \bfx,
\\
\M{C}^{\theta w}_{n,e} & = & \int_{T_e} h_\alpha(\theta_{n,e}^h( \bfx ))
\M{N}_e( \bfx )\trn \M{N}_e( \bfx ) \de \bfx,
\\
\M{C}^{\theta\theta}_{n,e} & = & \int_{T_e} \left( \rho_s C_s(
\theta_{n,e}^h( \bfx ) ) - h_d  w_d
 (\theta_{n,e}^h( \bfx )) \right) \M{N}_e( \bfx )\trn \M{N}_e(
\bfx ) \de \bfx,
\end{eqnarray}
\end{subequations}
whereas the symmetric positive-semidefinite blocks
$\M{K}_{\bullet}^{\bullet} \in \mathbb R^{3\times 3}$ attain the
form
\begin{subequations}
\begin{eqnarray}
\M{K}^{wP}_{n,e} & = & \int_{T_e} \frac{1}{g} \kappa \left(
  \theta_{n,e}^h( \bfx ),
  P_{n,e}^h( \bfx )
\right) \left( \nabla \M{N}_e( \bfx ) \right) \trn \nabla \M{N}_e(
\bfx ) \de \bfx
\nonumber \\
& + & \int_{\partial T_e \cap ( \Gamma_N^h \cup \Gamma_R^h ) }
\beta_c\, \M{N}_e(\bfx) \trn \M{N}_e(\bfx) \de \bfS,
\\
\M{K}^{\theta\theta}_{n,e}
& = &
\int_{T_e} \frac{1}{g} \lambda_c
\left(
  \theta_{n,e}^h( \bfx ),
  P_{n,e}^h( \bfx )
\right) \left( \nabla \M{N}_e (\bfx) \right) \trn \nabla
\M{N}_e(\bfx) \de \bfx
\nonumber \\
& + & \int_{\partial T_e \cap ( \Gamma_N^h \cup \Gamma_R^h ) }
\alpha_c \, \M{N}_e(\bfx) \trn \M{N}_e(\bfx) \de \bfS.
\end{eqnarray}
\end{subequations}
The non-linear terms $\M{R}_e^\bullet \in \mathbb{R}^{3 \times 1}$ are
expressed as
\begin{subequations}
\begin{eqnarray}
\M{R}_e^w( \M{\theta} )
& = &
\frac{1}{\Delta t}
\int_{T_e}
w_d \left( \M{N}_e(\bfx) \M{\theta} \right)
\M{N}_e(\bfx)\trn
\de \bfx,
\\
\M{R}_e^\theta( \M{\theta} )
& = &
\int_{\partial T_e \cap \Gamma_R^h}
e \sigma \left( \M{N}_e(\bfx) \M{\theta} \right)^4
\M{N}_e(\bfx) \trn \de \bfx
\end{eqnarray}
\end{subequations}
and the right hand-side-blocks $\M{F}_\bullet^\bullet \in \mathbb
R^{3\times 1}$ are provided by
\begin{subequations}
\begin{eqnarray}
\M{F}_{n+1,e}^w
& = &
\frac{1}{\Delta t}
\int_{T_e}
w_d \left( \theta_{n,e}^h( \bfx ) \right)
\M{N}_e(\bfx) \trn
\de \bfx
\nonumber \\
& + &
\int_{\partial T_e \cap ( \Gamma_N^h \cup \Gamma_R^h ) }
\beta_c P_\infty( t_{n+1} )
\M{N}_e(\bfx) \trn
\de \bfS,
\\
\M{F}_{n+1,e}^\theta & = & \int_{T_e} \frac{C_w}{g} \kappa\left(
\theta_{n,e}^h( \bfx ), P_{n,e}^h( \bfx ) \right)
\nabla\theta_{n,e}^h( \bfx ) \cdot \nabla P_{n,e}^h( \bfx )
\M{N}_e(\bfx)\trn \de \bfx
\nonumber \\
& + &
\int_{\partial T_e \cap ( \Gamma_N^h \cup \Gamma_R^h ) }
\alpha_c \theta_\infty( t_{n+1} ) \M{N}_e(\bfx) \trn
\de \bfS
\nonumber \\
& + & \int_{\partial T_e \cap \Gamma_R^h}
e\sigma
\theta_\infty^4( t_{n+1}) \M{N}_e(\bfx)\trn \de \bfS.
\end{eqnarray}
\end{subequations}
In the following example, the integrals were approximated
using the~$5\times5$-point Gauss quadrature and the non-linear
system~\eqref{eq:discrete_system} was solved iteratively using the
Newton method with the residual tolerance set to $10^{-8}$.

\begin{figure}[ht]
\begin{minipage}[b]{.47\textwidth}
\includegraphics*[width=5.2cm]{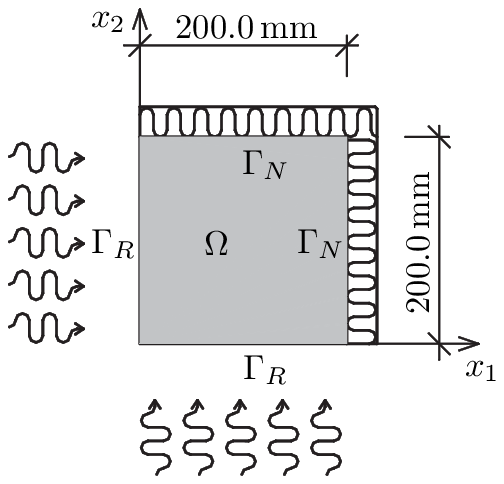}
\caption{Cross-section of specimen}\label{concrete_column}
\end{minipage}
\hfill
\begin{minipage}[b]{.52\textwidth}
\includegraphics[angle=270,width=\textwidth]{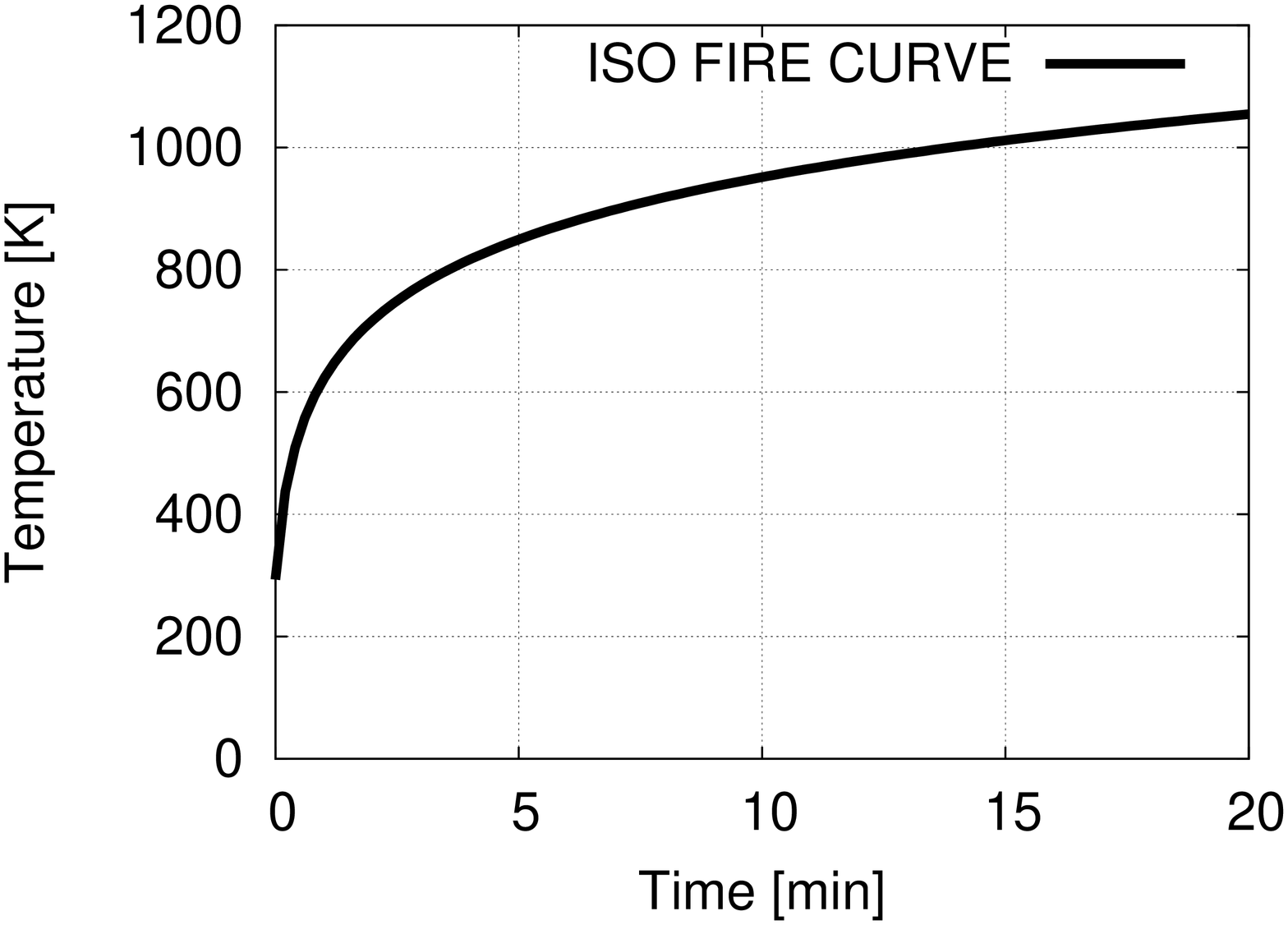}
\caption{ISO fire curve}\label{ISO fire curve}
\end{minipage}
\end{figure}

\subsection{Example} \label{example}
Our model problem deals with a square concrete specimen of
cross-section $200\times200~{\rm mm^2}$, which is exposed to fire on
the part $\Gamma_R$ of the boundary and insulated on the remaining
portion $\Gamma_N$ (see Fig.~\ref{concrete_column}). We assume that
the ambient temperature in the vicinity of $\Gamma_R$ increases
according to the ISO fire curve,
${\theta_{\infty}(t)=345\log(8t+1)+298.15}$, with $t$ given in
minutes~(see Fig.~\ref{ISO fire curve}).

The uniform initial conditions are set to
$\theta_0=298.15~{\rm K}$ and $P_0=2.7542\times10^3~{\rm Pa}$ and the constants
appearing in boundary conditions on $\Gamma_R$ are taken as
$\alpha_c=25~{\rm{W\,m^{-2}K^{-1}}}$, $\beta_c=0.019~{\rm{m\,s^{-1}}}$,
$P_{\infty}(t)=2.7542\times10^3~{\rm Pa}$,
$\sigma=5.67\times10^{-8}~{\rm W\,m^{-2}K^{-4}}$ and $e=0.7$.

\subsubsection{Material data for concrete at high
temperatures}\label{Material_data}

\begin{table}[t]
\centering
\begin{tabular}{llrlcc}
\hline
Physical quantity & Notation & Value & Dimension
\\
\hline\hline
Density of the solid microstructure  & $\rho_s$ & 2400.0 & ${\rm
kg\, m^{-3}}$ \\
Specific heat of liquid phase & $C_w$ & 4181.0 & ${\rm J\, kg^{-1}
K^{-1}}$ \\
Enthalpy of dehydration & $h_d$ & $2.4 \times 10^6$ & ${\rm J \, kg^{-1}}$
\\
Mass of cement per m$^3$ of concrete & $c$ & 300.0 & ${\rm kg\,
m^{-3}}$
\\
\hline
\end{tabular}
\caption{Material constants of concrete}\label{material_constants}
\end{table}

Basic material constants for concrete employed in this example
appear in Table~\ref{material_constants}. Following \cite{eurocode},
the thermal conductivity of concrete $\lambda_c$ can be bounded by
lower $\lambda_l$ and upper $\lambda_u$ limit values defined by
\begin{eqnarray*}
\lambda_l(\theta)&=&2.0-0.2451((\theta-273.15)/100.0)+0.0107((\theta-273.15)/100.0)^2,\\
\lambda_u(\theta)&=&1.36-0.136((\theta-273.15)/100.0)+0.0057((\theta-273.15)/100.0)^2
\end{eqnarray*}
and is set to $\lambda_c(\theta) =
(\lambda_l(\theta)+\lambda_u(\theta))/2$ in the results reported
below.

\begin{figure}[h]
\begin{minipage}[b]{.49\textwidth}
\includegraphics[angle=270,width=6.5cm]{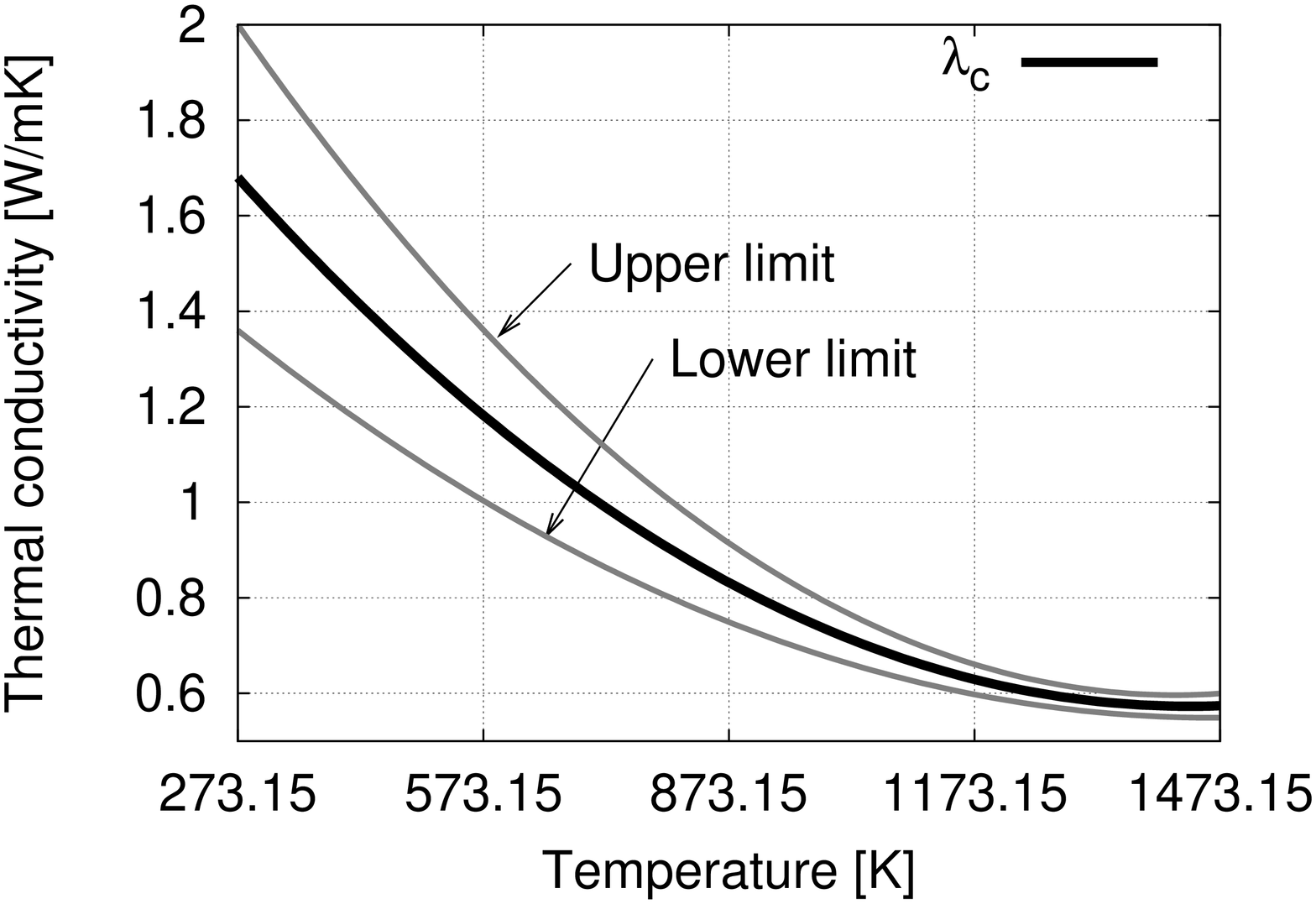}
\caption{Thermal conductivity of concrete}
\end{minipage}
\begin{minipage}[b]{.45\textwidth}
\includegraphics[angle=270,width=6.0cm]{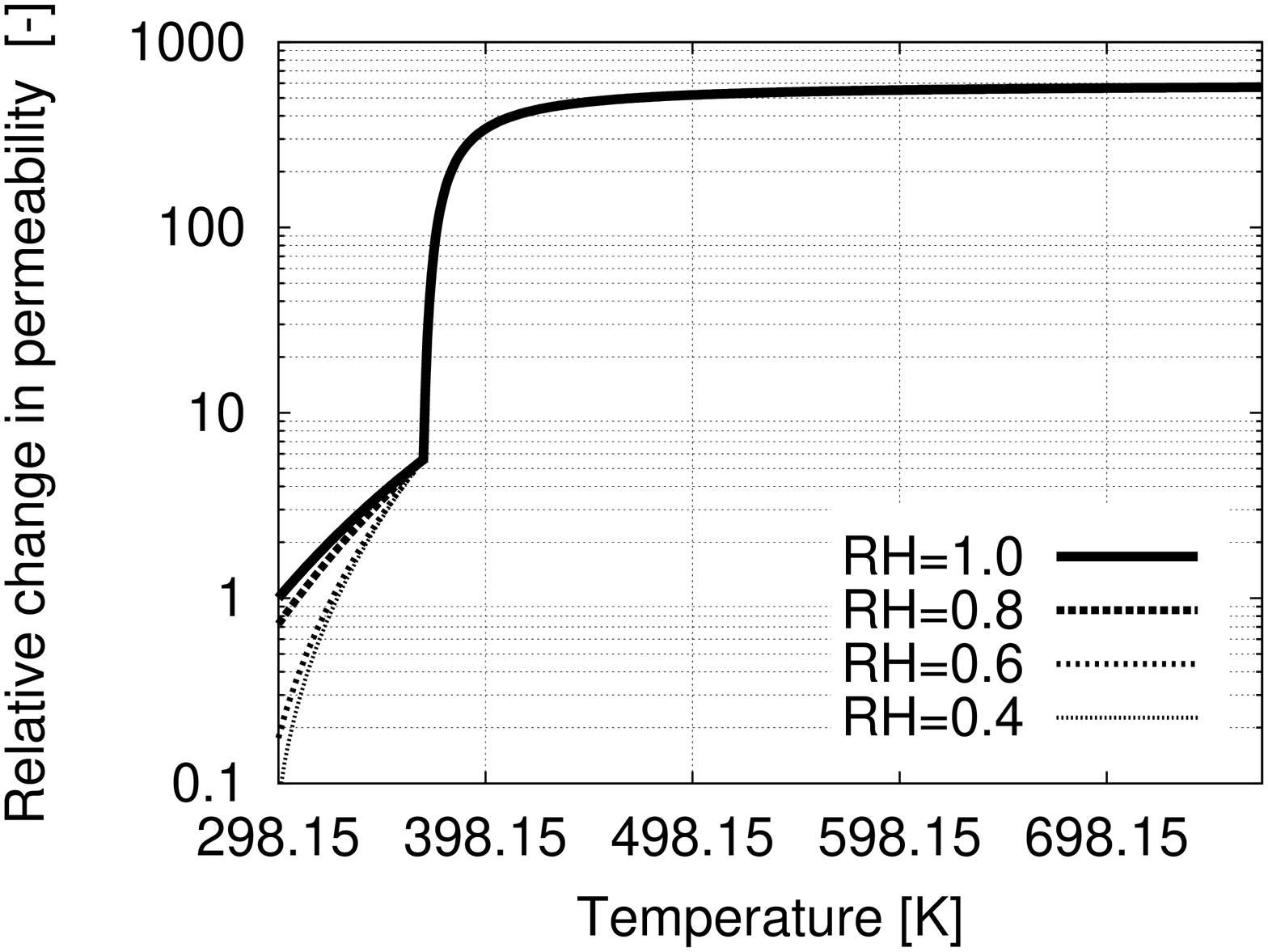}
\caption{Permeability of concrete}\label{permeability}
\end{minipage}
\end{figure}
\begin{figure}[h]
\begin{minipage}[b]{.49\textwidth}
\includegraphics[angle=270,width=6.3cm]{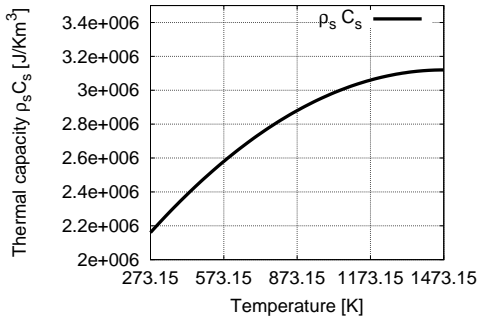}
\caption{Thermal capacity of solid skeleton}\label{capacity_solid}
\end{minipage}
\begin{minipage}[b]{.45\textwidth}
\includegraphics[angle=270,width=6.0cm]{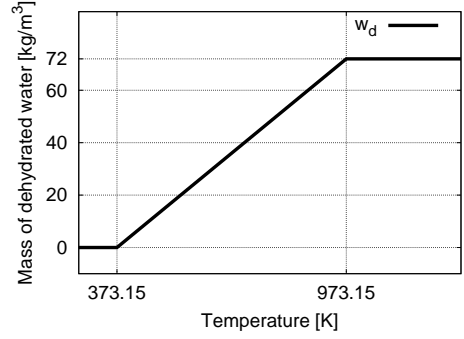}
\caption{Mass of dehydrated water}\label{dehydrated_water}
\end{minipage}
\end{figure}
The permeability of concrete $\kappa=\kappa(\theta,P)$ is adopted
from~\cite[(12a)--(12b)]{BaTh1978} and the relative change in
permeability $\kappa/\kappa_0$~[-], $\kappa_0=10^{-13}$~[ms$^{-1}$],
is displayed in Fig.~\ref{permeability} as a function of temperature
$\theta$ and relative humidity $RH$, defined as
$RH(P,\theta)=P/P_{sat}(\theta)$, where $P_{sat}(\theta)$ is the
vapor saturation pressure~(see~\cite{DaPeBi2006}).\footnote{\newtext{Note that
the pore pressure $P$ can generally exceed the saturation pressure $P_{sat}$, so that
$RH > 1$.}}

 The specific heat of solid matrix $C_s$ is
considered in the form (cf.~\cite{DaPeBi2006})
\begin{equation}
C_s(\theta)=900.0+80.0(\theta-273.15)/120.0-4.0((\theta-273.15)/120.0)^2.
\end{equation}
The thermal capacity of solid skeleton is displayed in Fig.
\ref{capacity_solid}. Assuming that concrete is fully hydrated at
room temperature, the mass of dehydrated water is given
as~\cite{Dwaikat:2009:HM}
\begin{equation*}
w_d( \theta ) = \left\{
 \begin{array}{rl}
 0.0 &
 \text{for } \theta \leq 373.15 \,\text{K},
 \\
0.04 c (\theta-373.15)/100.0 &
 \text{for } 373.15  < \theta \leq 973.15 \,\text{K},
 \\
 0.24 c &
 \text{for } \theta > 973.15 \,\text{K},
 \end{array}
\right.
\end{equation*}
see also Fig.~\ref{dehydrated_water} for an illustration. The
temperature dependence of the enthalpy of evaporation $h_{\alpha}$
will be approximated by the Watson equation~\cite{GaMaSch1999}
\begin{equation}\label{watson}
h_{\alpha}(\theta)=2.672\times10^5(\newtext{647.3-\theta})^{0.38}
\end{equation}
provided $\theta \leq 647.3$ K. Note that for higher temperatures
there is no liquid water in the pores and $h_{\alpha}(\theta)=0$.

The last comment concerns the sorption isotherms. For relative
humidities $RH<0.96$ and $RH>1.04$, thermodynamics-based relations
$w=\Phi(\theta,P)$ introduced in~\cite{BaTh1978} are adopted. In the
transition range, we employ the $C^1$-continuous cubic
interpolation.

\begin{figure}[p]
\begin{tabular}{cc}
\includegraphics*[width=6.5cm]{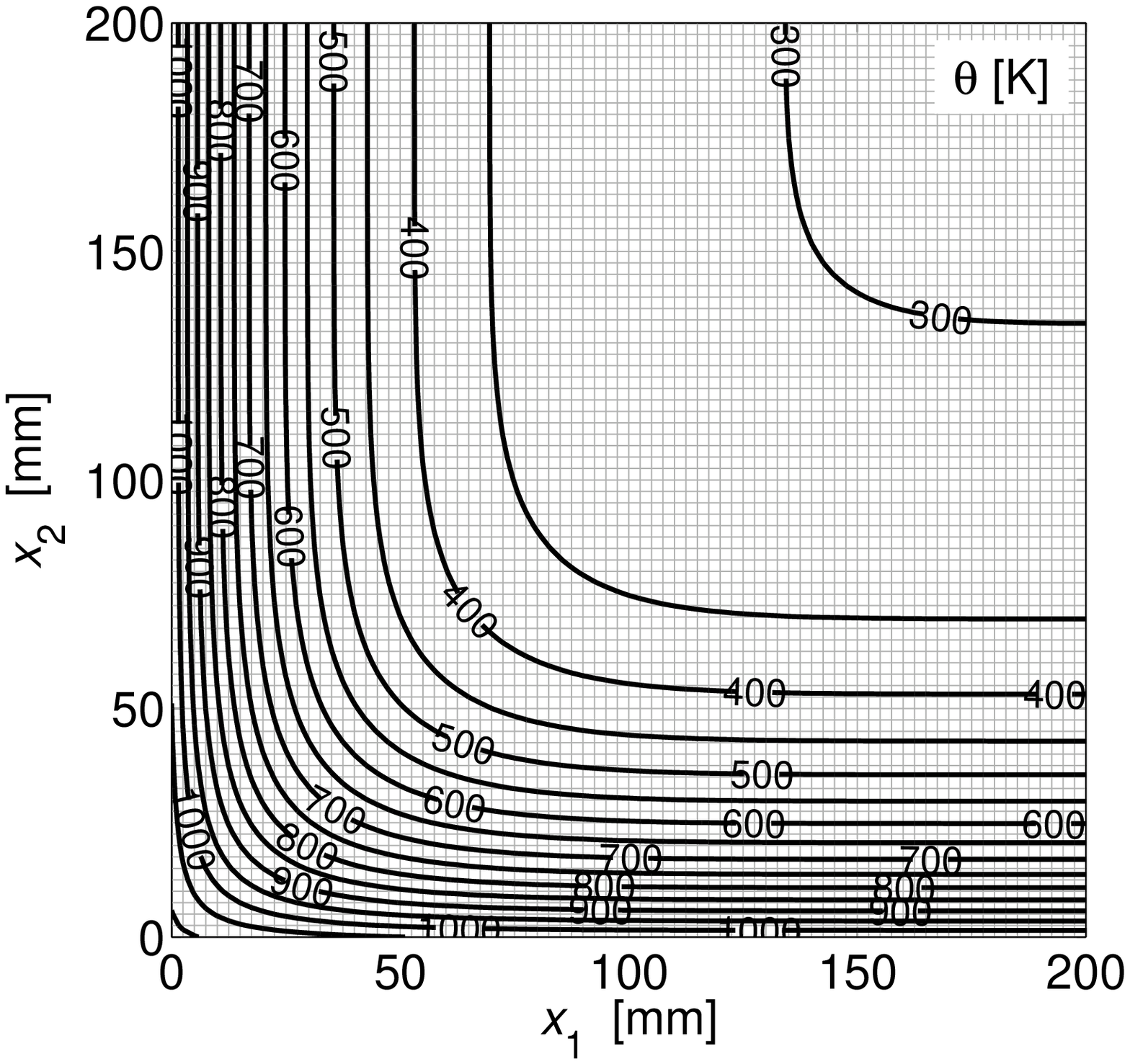}
&
\includegraphics*[width=6.5cm]{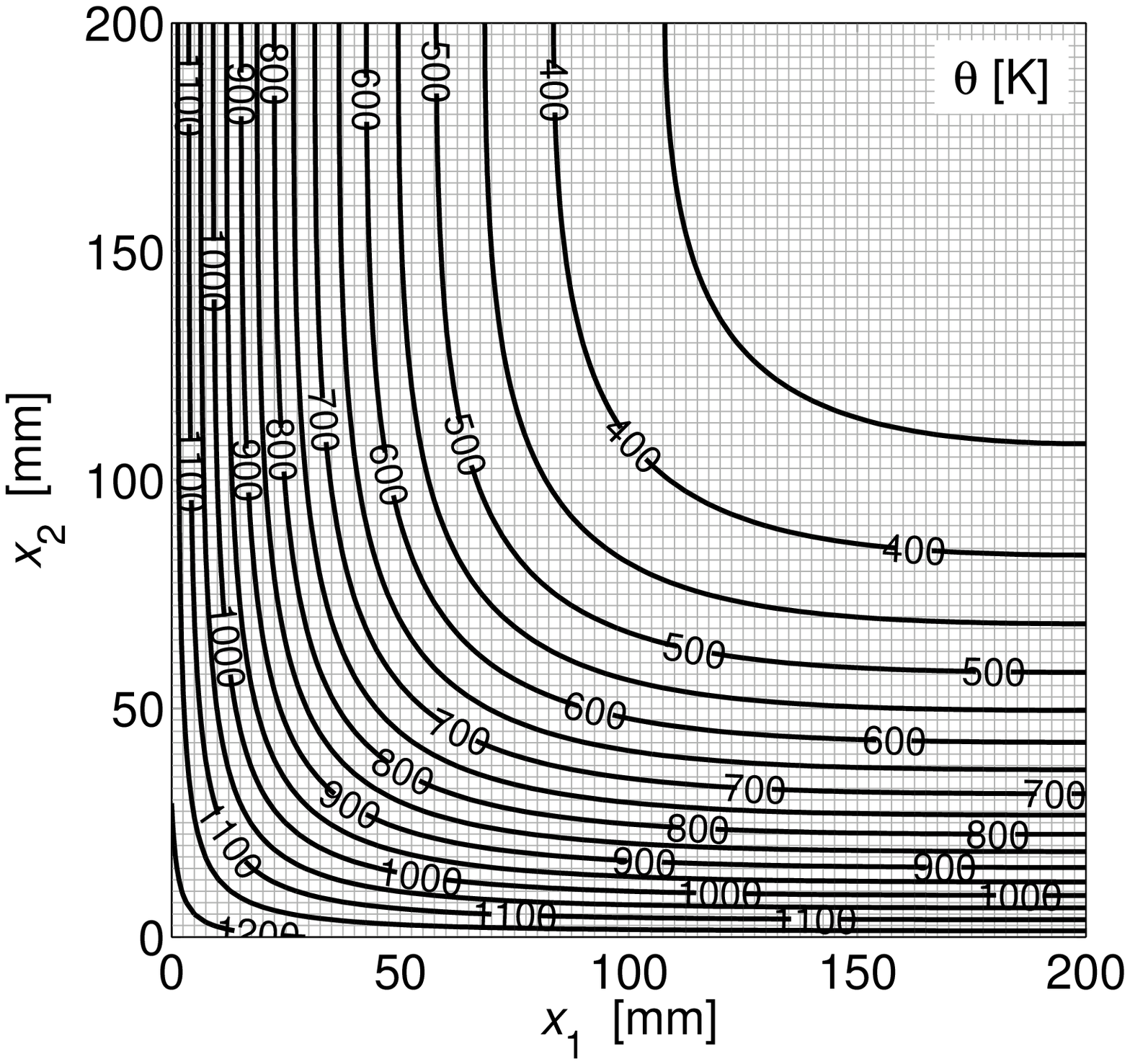}
\\
\includegraphics*[width=6.5cm]{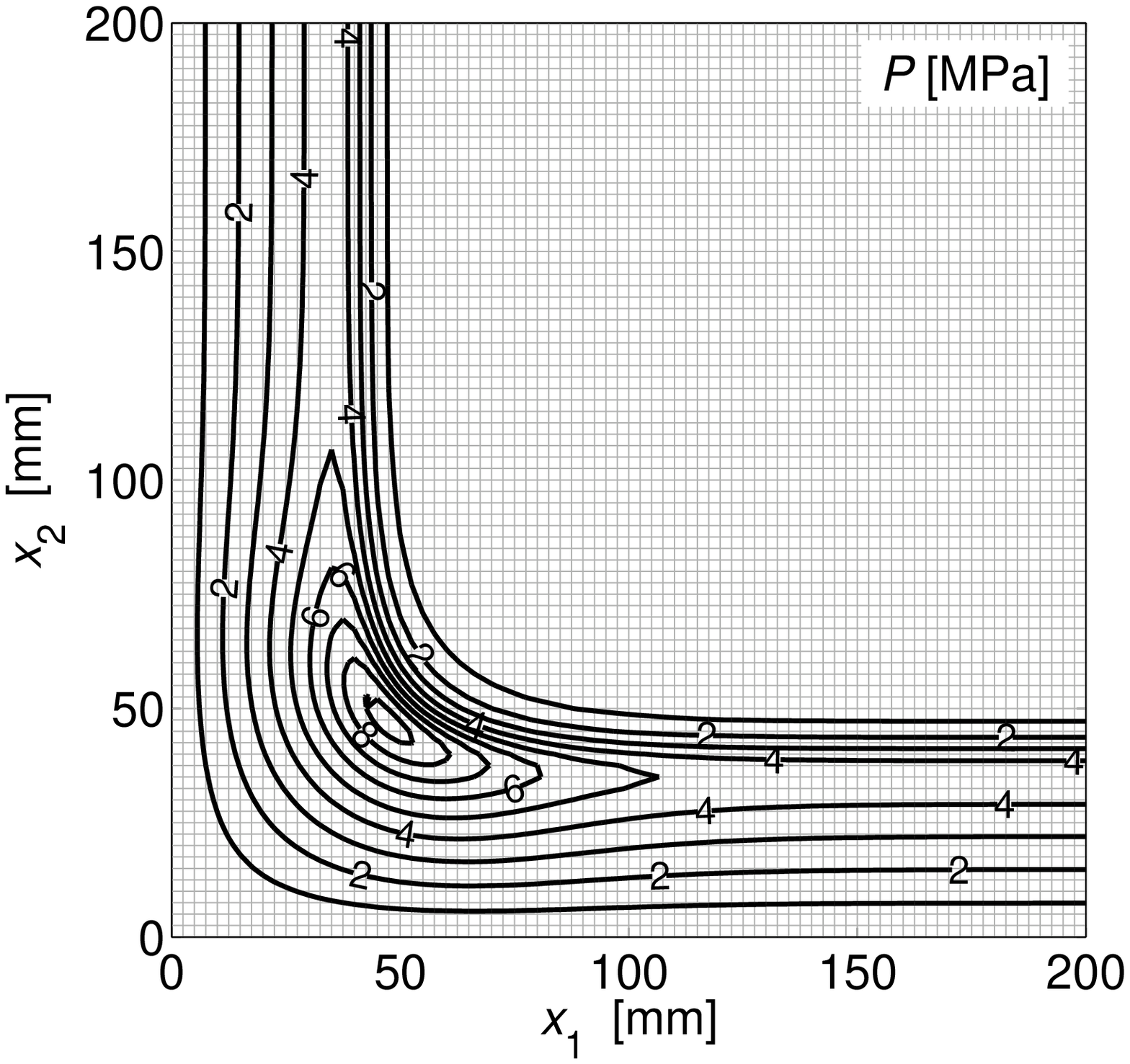}
&
\includegraphics*[width=6.5cm]{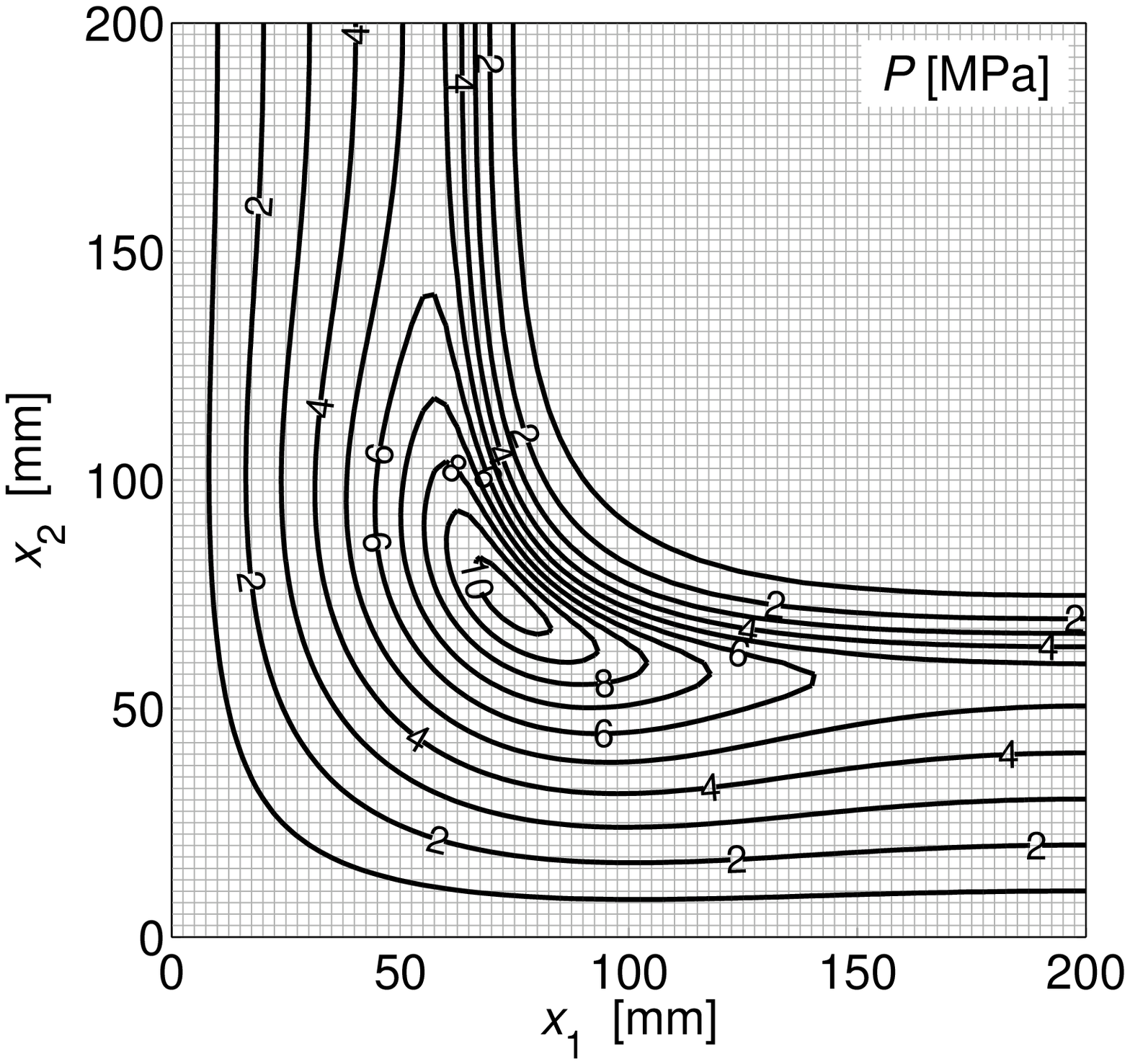}
\\
\includegraphics*[width=6.5cm]{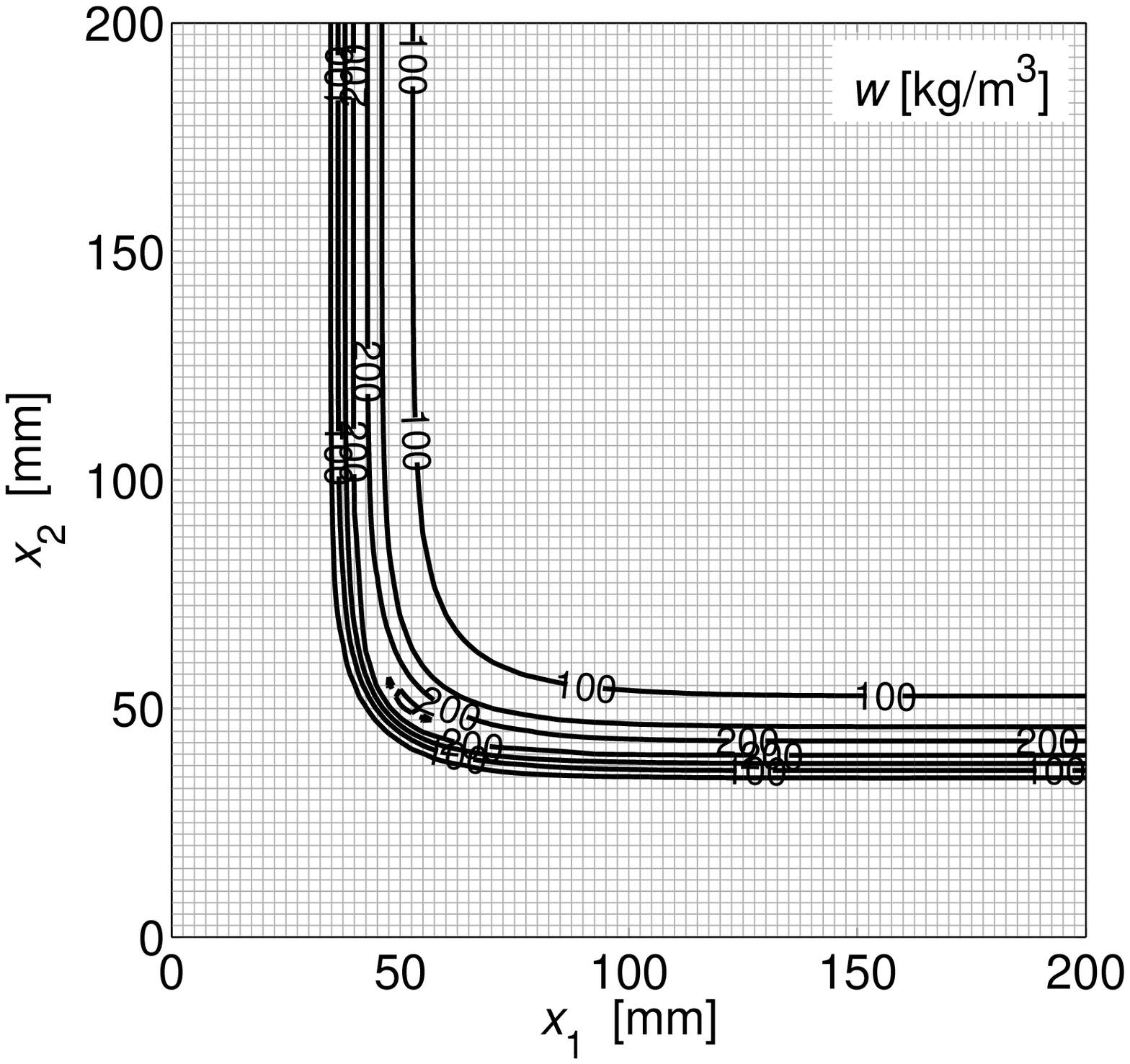}
&
\includegraphics*[width=6.5cm]{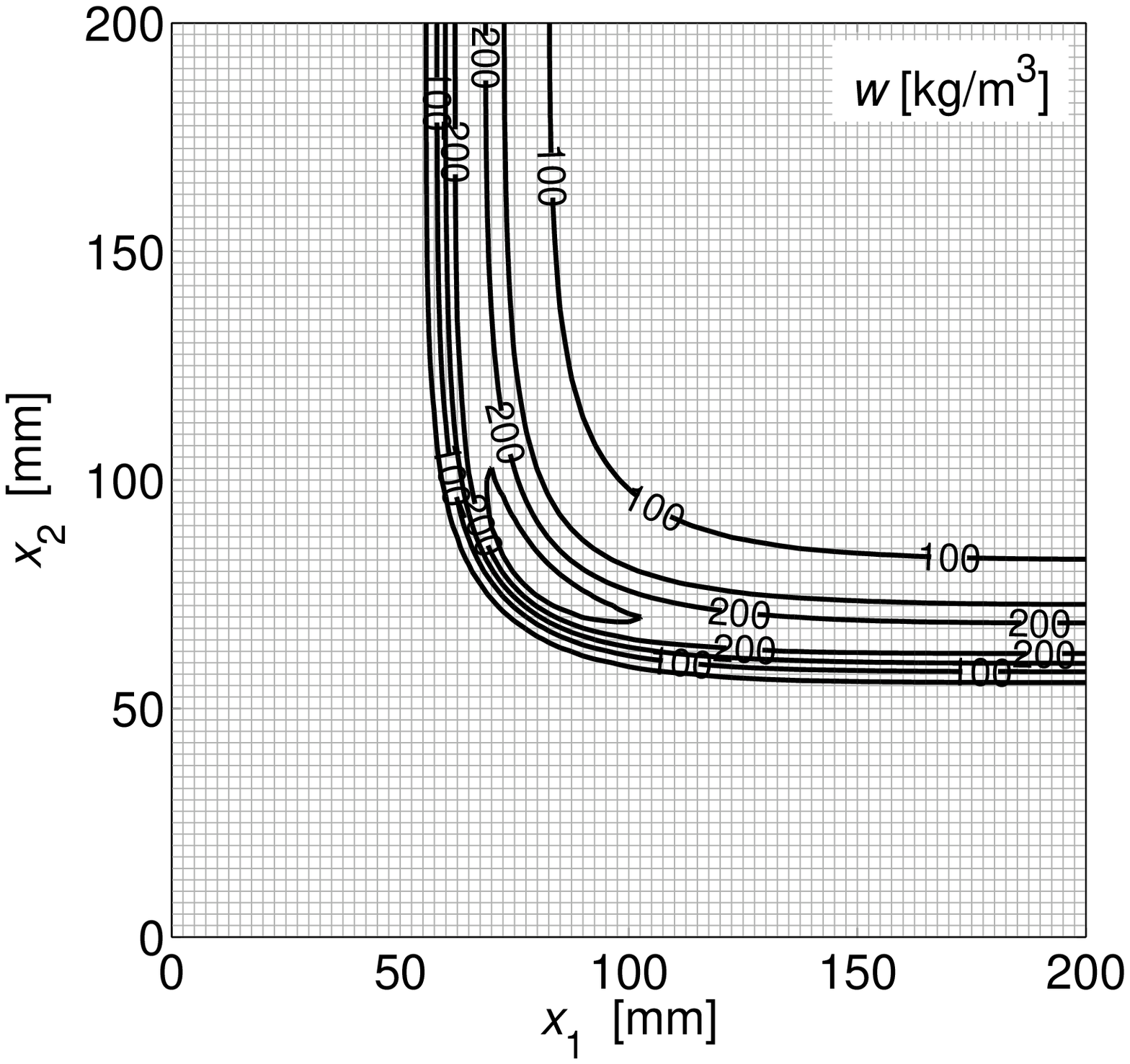}
\end{tabular}
\caption{Temperature, pore pressure and water content distribution
30 min (left) and 60 min (right) after exposure to fire.}
\label{fig_temp_press_water}
\end{figure}

\subsubsection{Results}
The results presented in this section are obtained using an in-house MATLAB code
and correspond to the \newtext{uniform} spatial discretization by $80\times80$
square elements~\newtext{(without any adaptivity)} and to the time step $\Delta
t=5$~s. The distribution of individual fields at two characteristic times
appears in Fig.~\ref{fig_temp_press_water}.

We observe that the numerical model correctly reproduces the rapid
heating of concrete specimen, accompanied by highly localized
profile of water content distribution. Mainly the latter phenomenon,
which is accurately resolved by the adopted fine grid, then
contributes to the development of high values of pore pressure in
the interior of the structure, leading to a potential threat to its
stability during fire due to explosive spalling.

\begin{figure}[h]
\begin{tabular}{cc}
\includegraphics*[width=6.5cm]{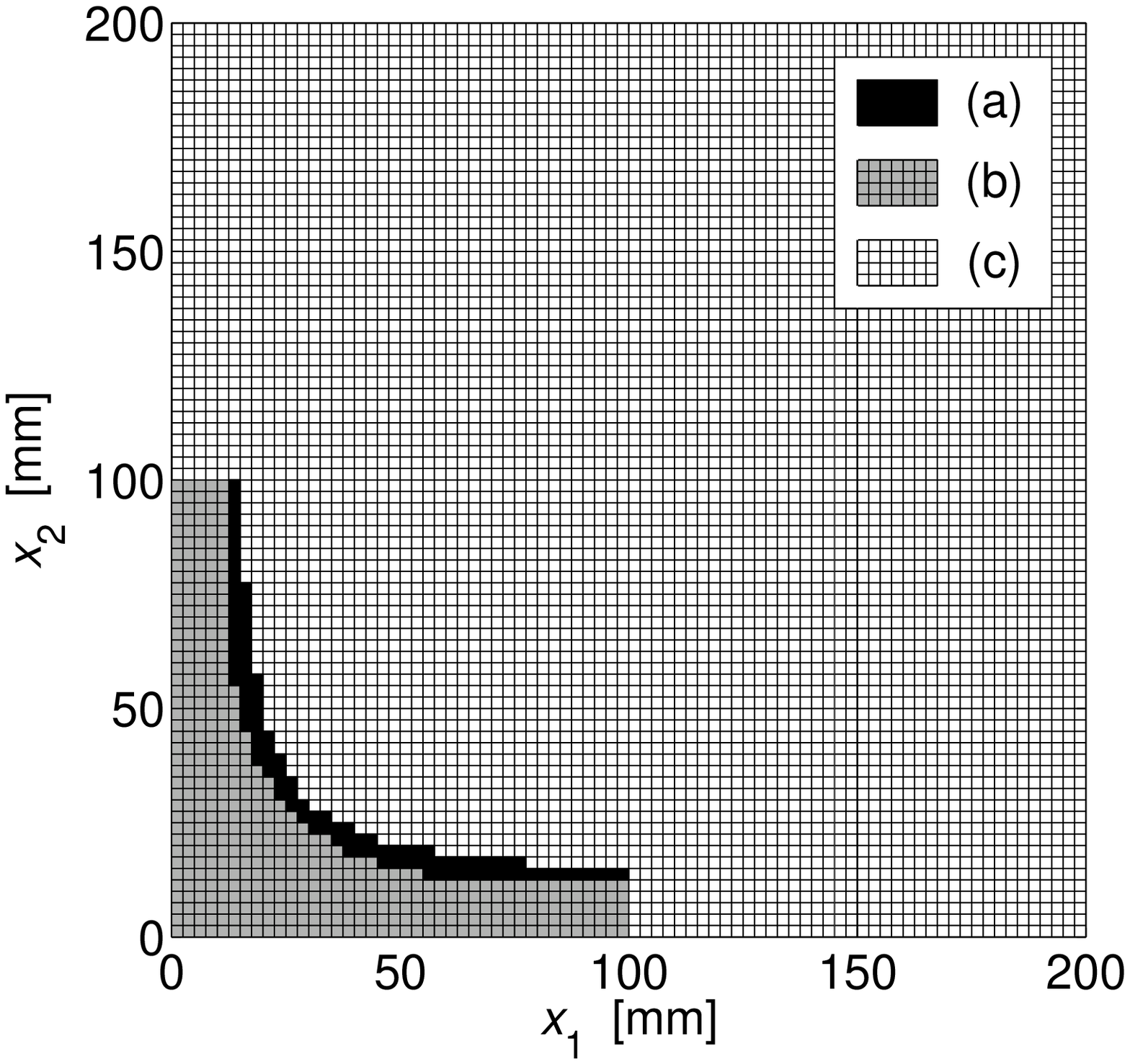}
&
\includegraphics*[width=6.5cm]{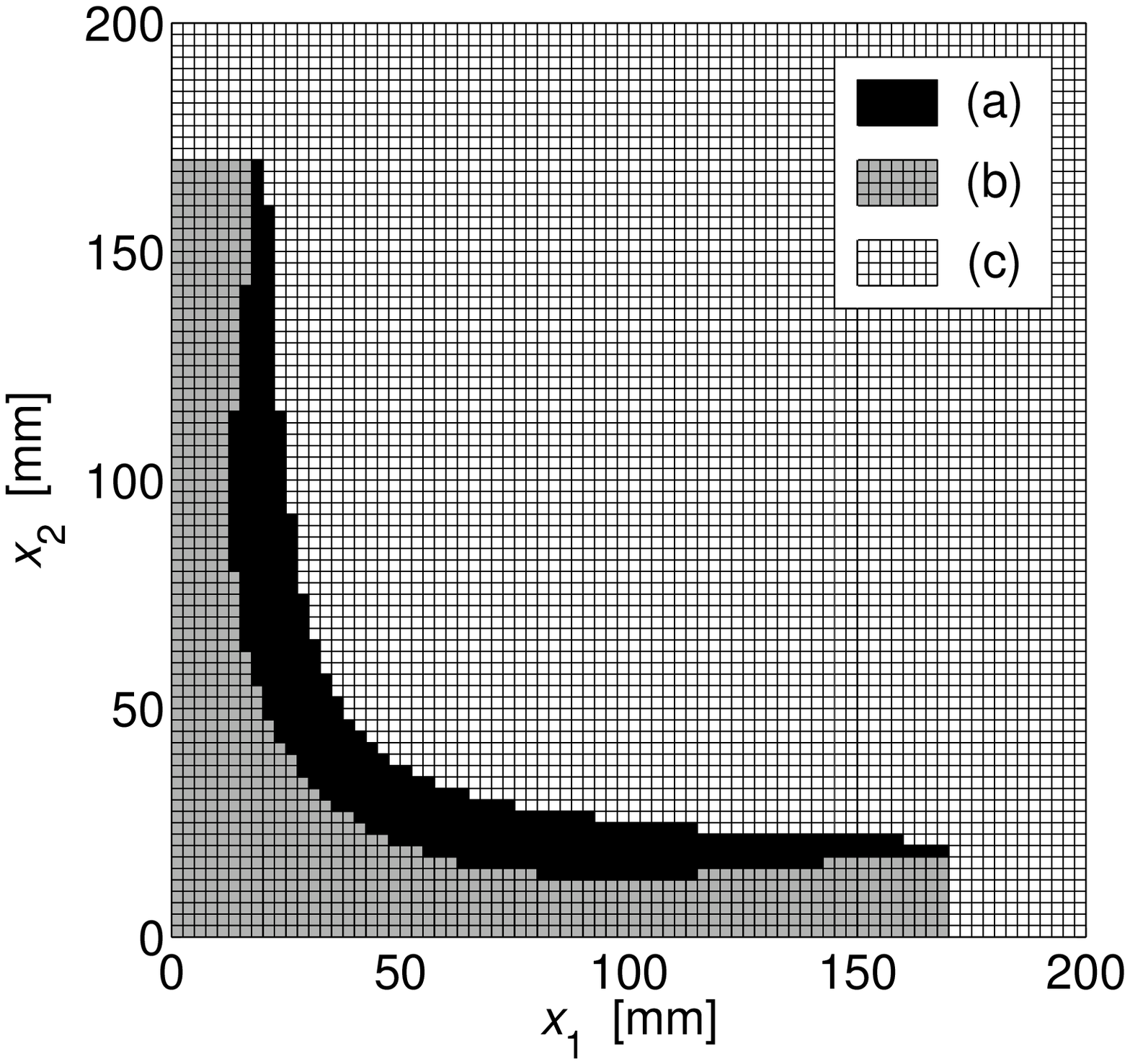}
\end{tabular}
\caption{Spalling of concrete 40 min (left) and 60 min (right) after
exposure to fire. \newtext{Region (a) denotes elements
failed due to spalling, (b)~marks unstable elements and (c)~stable part of the
cross-section.}}\label{fig_spalling}
\end{figure}

In order to assess this behavior more quantitatively, we adopt a
conservative engineering approach and assume that the spalling of
concrete occurs when the (appropriately reduced) pore pressure
exceeds the temperature-dependent tensile strength of concrete. In
particular, the spalling at $\bfx \in \Omega^h$ and time $t_n$
occurs when~\cite[Section~3.3]{Dwaikat:2009:HM}\footnote{\newtext{It should be
emphasized that relation~\eqref{eq:spalling_risk} is purely heuristic;
it nevertheless corresponds surprisingly well with detailed numerical
simulations cf.~\cite[Section 4.2]{Ozbolt2008}.}}
\begin{equation}\label{eq:spalling_risk}
{\phi} P^h_n ( \bfx ) \geq f_t \left( \theta^h_n(\bfx) \right),
\end{equation}
where $\phi$ is the porosity of concrete, \newtext{$P^h_n$ and $\theta^h_n$
denote the finite element approximations of pore pressure and temperature for
mesh size $h$ and the $n$-th time step} and $f_t$ as a function of $\theta$
is provided by a piecewise linear relation~\cite[Section~4.4]{Dwaikat:2009:HM}
\begin{equation}
f_t( \theta ) = f_{t0} \times \left\{
 \begin{array}{rl}
 1 &
 \text{for } \theta \leq 373.15~\text{K},
 \\
 (873.15 - \theta)/500 &
 \text{for } 373.15~\text{K} < \theta \leq 823.15~\text{K},
 \\
 (1473.15 - \theta)/6500 &
 \text{for } 823.15~\text{K} < \theta \leq 1473.15~\text{K},
 \\
 0 & \text{otherwise,}
 \end{array}
\right.
\end{equation}
where $f_{t0}$ designates the initial tensile strength of concrete.

The spatial distribution of the spalling damage for $\phi = 0.1$ and
$f_{t0} = 2$~MPa appears in Fig.~\ref{fig_spalling}. Here, three
different zones can be distinguished. The first region (a),
highlighted in black color, corresponds to elements failed due to
spalling damage as predicted by the
criterion~\eqref{eq:spalling_risk} for the element center. The
second zone (b) corresponds to the part of the structure in which
the local strength is sufficiently high to sustain the pore
pressure, but its stability is lost due the explosive spalling of
the former region. Finally, the light gray zone (c) indicates the
portion of the cross-section still capable of transmitting stresses
due to mechanical loading, which is thereby responsible for the
structural safety during fire.

\section*{Acknowledgments}
This outcome has been achieved with the financial support of the
Ministry of Education, Youth and Sports of the Czech Republic,
project No.~1M0579, within activities of the CIDEAS research centre
(the first author) and project No. MSM6840770001 (the second
author). Additional support from the grants 201/09/1544~(the first author), 
103/08/1531 and 201/10/0357 (the third author) provided by the Czech Science
Foundation is greatly acknowledged.

\end{document}